\documentclass[12pt]{amsart}
\usepackage{amsbsy}
\usepackage{graphicx,epsfig,subfigure,psfrag}
\begin{document}

\newtheorem{theorem}{Theorem}
\newtheorem{proposition}{Proposition}
\newtheorem{lemma}{Lemma}
\newtheorem{corollary}{Corollary}
\newtheorem{definition}{Definition}
\newtheorem{remark}{Remark}
\newcommand{\be}{\begin{equation}}
\newcommand{\ee}{\end{equation}}
\newcommand{\tex}{\textstyle}
\numberwithin{equation}{section} \numberwithin{theorem}{section}
\numberwithin{proposition}{section} \numberwithin{lemma}{section}
\numberwithin{corollary}{section}
\numberwithin{definition}{section} \numberwithin{remark}{section}
\newcommand{\ren}{\mathbb{R}^N}
\newcommand{\re}{\mathbb{R}}
\newcommand{\n}{\nabla}
\newcommand{\iy}{\infty}
\newcommand{\pa}{\partial}
\newcommand{\fp}{\noindent}
\newcommand{\ms}{\medskip\vskip-.1cm}
\newcommand{\mpb}{\medskip}
\newcommand{\BB}{{\bf B}}
\newcommand{\AAA}{{\bf A}}
\newcommand{\Am}{{\bf A}_{2m}}
\newcommand{\ef}{\eqref}
\renewcommand{\a}{\alpha}
\renewcommand{\b}{\beta}
\newcommand{\g}{\gamma}
\newcommand{\G}{\Gamma}
\renewcommand{\d}{\delta}
\newcommand{\D}{\Delta}
\newcommand{\e}{\varepsilon}
\newcommand{\var}{\varphi}
\renewcommand{\l}{\lambda}
\renewcommand{\o}{\omega}
\renewcommand{\O}{\Omega}
\newcommand{\s}{\sigma}
\renewcommand{\t}{\tau}
\renewcommand{\th}{\theta}
\newcommand{\z}{\zeta}
\newcommand{\wx}{\widetilde x}
\newcommand{\wt}{\widetilde t}
\newcommand{\noi}{\noindent}
\newcommand{\uu}{{\bf u}}
\newcommand{\UU}{{\bf U}}
\newcommand{\VV}{{\bf V}}
\newcommand{\ww}{{\bf w}}
\newcommand{\vv}{{\bf v}}
\newcommand{\WW}{{\bf W}}
\newcommand{\hh}{{\bf h}}
\newcommand{\di}{{\rm div}\,}
\newcommand{\inA}{\quad \mbox{in} \quad \ren \times \re_+}
\newcommand{\inB}{\quad \mbox{in} \quad}
\newcommand{\inC}{\quad \mbox{in} \quad \re \times \re_+}
\newcommand{\inD}{\quad \mbox{in} \quad \re}
\newcommand{\forA}{\quad \mbox{for} \quad}
\newcommand{\whereA}{,\quad \mbox{where} \quad}
\newcommand{\asA}{\quad \mbox{as} \quad}
\newcommand{\andA}{\quad \mbox{and} \quad}
\newcommand{\withA}{,\quad \mbox{with} \quad}
\newcommand{\orA}{,\quad \mbox{or} \quad}
\newcommand{\ssk}{\smallskip}
\newcommand{\LongA}{\quad \Longrightarrow \quad}
\def\com#1{\fbox{\parbox{6in}{\texttt{#1}}}}
\def\N{{\mathbb N}}
\def\A{{\cal A}}
\def\WW{{\cal W}}
\newcommand{\de}{\,d}
\newcommand{\eps}{\varepsilon}
\newcommand{\spt}{{\mbox spt}}
\newcommand{\ind}{{\mbox ind}}
\newcommand{\supp}{{\mbox supp}}
\newcommand{\dip}{\displaystyle}
\newcommand{\prt}{\partial}
\renewcommand{\theequation}{\thesection.\arabic{equation}}
\renewcommand{\baselinestretch}{1.2}

\title
{\bf Classification of global and  blow-up\\  sign-changing
solutions of a semilinear heat \\ equation in the subcritical
Fujita
 range II. Higher-order diffusion}

\author{
V.A.~Galaktionov, E.~Mitidieri,  and S.I.~Pohozaev}

\address{Department of Mathematical Sciences, University of Bath,
 Bath BA2 7AY, UK}
\email{vag@maths.bath.ac.uk}

\address{Dipartimento di Matematica e Informatica,
Universit\`a  di Trieste, Via  Valerio 12, 34127 Trieste, ITALY}
 \email{mitidier@units.it}

\address{Steklov Mathematical Institute,
 Gubkina St. 8, 119991 Moscow, RUSSIA}
\email{pokhozhaev@mi.ras.ru}



 \date{\today}




\begin{abstract}

A detailed study of two classes of oscillatory global (and
blow-up) solutions was began in \cite{GMPOscI} for
 the semilinear heat equation in the subcritical Fujita range
  \be
   \label{0111A}
    \tex{
u_t = \D u +|u|^{p-1}u \quad \mbox{in} \quad \ren \times \re_+
\quad \mbox{for}
 \,\,\, 1 < p \le p_0=1+ \frac 2N,
 }
  \ee
  with bounded integrable initial data $ u(x,0)=u_0(x)$. This
  study
 is continued and extended here
for the $2m$th-order heat equation, for   $m \ge 2$, with
non-monotone nonlinearity
  \be
   \label{111}
    \tex{
u_t = -(-\D)^m u +|u|^p \quad \mbox{in} \quad \ren \times \re_+,
\quad \mbox{in the range}
 \,\,\, 1 < p \le p_0=1+ \frac {2m}N,
 }
 \ee
with the same initial data $u_0$. The fourth order bi-harmonic
case $m=2$ is studied in greater detail. The blow-up Fujita-type
result for \ef{111} now reads as follows:
  blow-up occurs for any initial data $u_0$ with  positive first
Fourier coefficient:
 $$
\tex{
  \int u_0(x)\,{\mathrm d} x>0,
  }
  $$
  i.e.,  as for \ef{0111A}, any such arbitrarily small initial function $u_0(x)$ leads to blow-up. The construction
  of two countable families of global sign changing solutions is performed
  on the basis of bifurcation/branching analysis and a further analytic-numerical study.
  In particular, a countable sequence of {\em bifurcation points} of similarity solutions is
  obtained:
  $$
 \mbox{$
p_l=1 + \frac {2m}{N+l}, \quad l=0,1,2, ... \, .
 $}
$$ 

\end{abstract}

\maketitle



\setcounter{equation}{0}
\section{Introduction: higher-order semilinear heat equations, blow-up, Fujita exponent,
 and global oscillatory solutions}
 \label{Sect1}
  \setcounter{equation}{0}







 This paper generalizes  to higher-order (poly-harmonic) diffusion operators the study of \cite{GMPOscI} of
  the semilinear second-order heat equation in the subcritical Fujita range
  \be
   \label{0111}
    \tex{
u_t = \D u +|u|^{p-1}u \quad \mbox{in} \quad \ren \times \re_+,
\quad \mbox{where}
 \,\,\, 1 < p \le p_0=1+ \frac 2N,
 }
  \ee
  with bounded integrable initial data $ u(x,0)=u_0(x)$.
 Necessary key references, results, and our main motivation of
 the study of \ef{0111} and further related models can be found in \cite[\S~1]{GMPOscI}.

\subsection{Higher-order semilinear heat equation:  blow-up  Fujita-like result}

Thus, we intend to extend some of key results of \cite{GMPOscI} to
semilinear parabolic equations with higher-order diffusion
operators. Such models are steadily becoming more and more popular
in various applications and in general PDE theory.
  Namely, we consider the
$2m$th-order heat equation,  for  $m \ge 2$, in the subcritical
Fujita range:
 \be
   \label{00111}
    \tex{
u_t = -(-\D)^m u +|u|^p \quad \mbox{in} \quad \ren \times \re_+
\quad \mbox{for}
 \,\,\, 1 < p \le p_0=1+ \frac {2m}N,
 }
  \ee
with sufficiently smooth, bounded, and integrable initial data,
\be
\label{in1}
 u(x,0)=u_0(x) \inB \ren.
 \ee
The choice of the non-monotone nonlinearity $|u|^p$ (a source
term) in \ef{00111} is associated with the necessity of having a
standard sounding Fujita-type blow-up result.
Namely, it is known \cite{EGKP1} that blow-up occurs for
\ef{00111} for any solutions with initial data having positive
first Fourier coefficient (see \cite{Eg4} for further details and
\cite{GPInd} for an alternative proof):
 \be
 \label{u01}
 \tex{
  \int_{\ren} u_0(x)\,{\mathrm d} x>0,
  }
  \ee
  i.e., any, even arbitrarily small, such data lead to blow-up. A similar result for positive solutions of
  \ef{0111} was well-known since Fujita work in 1966; see
  \cite[\S~1]{GMPOscI} for a survey.

\subsection{Results and layout of the paper}

  In
 Sections \ref{S2}--\ref{SectN4},
 we perform
   construction
  of countable sets of global sign changing solutions
  on the basis of bifurcation/branching analysis, as well as of a centre-stable manifold one.
Here, we apply spectral theory of related non-self-adjoint
$2m$th-order operators  \cite{Eg4}, which is available for any
$m=2,3,...$\,. As for \ef{0111}, i.e., for $m=1$, this gives a
similar sequence of critical bifurcation exponents:
 \be
 \label{h2}
 \mbox{$
p_l=1 + \frac {2m}{N+l}, \quad l=0,1,2, ... \, .
 $}
 \ee

 References and some results for  analogous  global similarity solutions
 of a different   higher-order reaction-diffusion PDE with a
 standard
 monotone nonlinearity as in \ef{0111},
 \be
 \label{2mGH}
 u_t= -(-\D)^m u + |u|^{p-1}u \quad (m \ge 2)
  \ee
 can be found in \cite{GHUni}.
 It is remarkable (and rather surprising for us)  that the bifurcation-branching phenomena therein for
 \ef{2mGH} are
 entirely different from those for  the present equation \ef{00111}, which turn out also to be more complicated,
 with various {\em standard} and {\em non-standard} bifurcation phenomena.

It is worth mentioning here that our study also directly  concerns
blow-up solutions of \ef{0111}: we claim that, under the
conditions that our two classes of its global oscillatory
solutions are {\em evolutionary complete} (see
\cite[\S~7]{GMPOscI} for a precise statement and some results for
\ef{0111}), all other solutions of \ef{00111} must blow-up in
finite time. Then this describes a much wider class of blow-up
solutions, and actually says that {\em almost all}\, (with {\em
a.a.} defined in a natural way) solutions of \ef{00111} in the
subcritical Fujita range blow-up in finite time.







\section{Global similarity solutions and $p$-bifurcation branches}
 \label{S2}



In what follows, we use a general scheme and ``ideology" of the
study in \cite{GMPOscI}
 of the second-order semilinear equation \ef{0111}.
 Therefore, omitting some obvious details,
 we now
more briefly start to describe which results on global solutions
can be extended to the $2m$th-order reaction-diffusion
equation \ef{00111} in the subcritical range.

\subsection{First family of global patterns: similarity solutions}

As usual, for the higher-order model (\ref{00111}) with $m \ge 2$,
we first study the existence
 and multiplicity of the standard global (i.e., well defined for all $t>0$)
  similarity solutions of the form
 \be
 \label{1.6R}
 u_S(x,t) = t^{-\frac 1{p-1}} f(y), \quad y = x/t^{\frac 1{2m}}.
 \ee
This leads to the semilinear elliptic problem for the rescaled
similarity profile $f$:
 \be
 \label{1.7R}
 \left\{
  \begin{matrix}
\BB_1 f +  |f|^p \equiv  -(-\D)^m f + \textstyle{\frac 1{2m}}\, y
\cdot \nabla f + \textstyle{\frac 1{p-1}} \, f + |f|^p=0 \quad
\mbox{in} \quad \ren,\ssk\ssk
\\
f(y) \,\,\, \mbox{decays exponentially fast as} \,\,\, |y| \to
\infty. \qquad\qquad\qquad\qquad\quad\,\,
\end{matrix}
\right.
 \ee

For $m=1$, this problem admits a variational setting in a weighted
metric of $L^2_\rho(\ren)$, where $\rho={\mathrm e}^{|y|^2/4}$.
This positive fact was heavily used in \cite{GMPOscI}, where
category/fibering techniques allowed us to detect a countable
number of solutions and bifurcation branches.


However, for any $m \ge 2$, \ef{1.7R} {\em is not variational} in
any weighted $L^2$ space; cf. reasons for that and a similar
negative result in \cite[\S~7]{GW2}. So, those power tools of
potential operator theory  in principle cannot be applied for
\ef{1.7R}, with any $m \ge 2$.

 Moreover,
unlike the previous study of \ef{0111} in \cite{GMPOscI}, we
cannot use standard variational results on bifurcation from
eigenvalues of arbitrary multiplicity (for our purposes, the
results for odd multiplicity \cite[p.~381, 401]{Deim} concerning
local and global continuation of branches are sufficient). We also
do not have global multiplicity results via
Lusternik--Schnirel'man (L--S, for short) and fibering theory. As
usual,  higher-order semilinear elliptic equations such as
\ef{1.7R}, or even the corresponding ODEs for radially symmetric
profiles $f$, become principally different and more difficult than
their second-order variational counterparts.
 We again refer to \cite[\S~6.7C]{Berger}
 for general results on bifurcation diagrams, and to \cite{Rynn2,
 Rynn1} for more detailed results for related $2m$th-order ODEs in
 1D. These results do not apply directly but can be used for a
 better understanding of global bifurcation diagrams of similarity patterns $f(y)$.

 Thus, as in \cite{GHUni} for a quite similar looking equation \ef{2mGH},
 for global continuation of branches,
  we have to rely more heavily  on
numerical methods, and this is an unavoidable feature of such a
study of nonlinear  higher-order equations. Surprisingly, we
detect {\em completely different} local and
 global properties of
$p$-branches in contrast with those in \cite{GHUni} for the
equation \ef{2mGH}, which therefore are not so definitely attached
to variational, monotone, or order-preserving (i.e., via the
Maximum Principle) features of these difficult global similarity
problems studied since the 1980s.


\subsection{Fundamental solution and Hermitian spectral theory}
 \label{Sect3}

 We begin with the necessary
  fundamental solution $b(x,t)$  of the corresponding linear parabolic (poly-harmonic) equation
  \be
  \label{Lineq}
  u_t = - (-\D)^m u \inA,
  \ee
 which takes the standard similarity form
  \be
  \label{1.3R}
   b(x,t) = t^{-\frac N{2m}}F(y) \whereA y= x/t^{\frac 1{2m}}.
 \ee
 The rescaled kernel $F$ is then the unique radial solution of the elliptic equation
  \begin{equation}
\label{ODEf}
 {\bf B} F \equiv -(-\Delta )^m F + \textstyle{\frac 1{2m}}\, y \cdot
\nabla F + \textstyle{\frac N{2m}} \,F = 0
 \,\,\,\, {\rm in} \,\, \ren,  \quad \mbox{with} \,\,\, \textstyle{\int F =
 1.}
\end{equation}
The rescaled kernel $F(|y|)$ is  oscillatory as $|y| \to \infty$
and satisfies the estimate
 \cite{EidSys, Fedor}
\begin{equation}
\label{fbar}
 |F(y)| < D\,\,  {\mathrm e}^{-d|y|^{\alpha}}
\,\,\,{\rm in} \,\,\, \ren, \quad \mbox{where} \,\,\,
\a=\textstyle{ \frac {2m}{2m-1}} \in (1,2),
\end{equation}
for some positive constants $D$ and $d$ depending on $m$ and $N$.
 The linear operator $\BB_1$ in equation
(\ref{1.7R}) is connected with the operator (\ref{ODEf}) for the
rescaled kernel $F$ in (\ref{1.3R}) by
 \be
 \label{BB*}
 \tex{
  \BB_1 = \BB + c_1 I, \quad \mbox{where} \,\,\, c_1 = \frac
  1{p-1}- \frac N{2m} \equiv
 \textstyle{\frac {N(p_0-p)}{2m(p-1)}} \andA p_0=1+\frac{2m}N.
 }
  \ee


In view of  (\ref{BB*}), in order to study the similarity
solutions,
 we need the  spectral properties
of ${\bf B}$ and of the corresponding adjoint operator ${\bf
B}^*$.
 Both are considered  in weighted $L^2$-spaces with the
weight functions induced by the exponential estimate of the
rescaled kernel (\ref{fbar}). For $m \ge 2$,
we consider $\BB$ in the weighted space $L^2_\rho(\ren)$ with the
exponentially growing weight function
 \be
  \label{rho44}
  \rho(y) = {\mathrm e}^{a |y|^\a}>0 \quad {\rm in} \,\,\, \ren,
 \ee
  where $a \in (0,  2d)$ is any fixed
constant and $d$ is as in \ef{fbar}. We ascribe  to ${\bf B}$ the
domain $H^{2m}_\rho(\ren)$ being
 a Hilbert space with the
norm $$
 \tex{
 \|v\|^2 = \int \rho(y)\sum\limits_{k=0}^{2m} |D^{k}
 v(y)|^2 \, \mathrm{d} y,
 }
  $$
induced by the corresponding inner product. Then $H^{2m}_{\rho}
\subset L^2_{\rho}  \subset L^2 $.  The spectral properties ${\bf
B}$
are as follows \cite{Eg4}:

\begin{lemma}
\label{lemspec} {\rm (i)} ${\mathbf B}: H^{2m}_\rho \to L^2_\rho$
is a bounded linear operator with the real point spectrum
\begin{equation}
\label{spec1} \sigma({\mathbf B}) = \{\lambda_l =
-\textstyle{\frac l{2m}}, \,\, l = 0,1,2,...\}.
\end{equation}
The eigenvalues $\l_l$ have finite multiplicities with
eigenfunctions
\begin{equation}
\label{eigen} \psi_\beta(y) =\textstyle{\frac{(-1)^{|\b|}}{\sqrt
{\b !}}} \, D^\beta F(y), \quad \mbox{for any} \,\,\,|\b|=l.
\end{equation}

{\rm (ii)} The set $\Phi = \{\psi_\b\}_{|\b| \ge 0}$ is complete
and the resolvent $({\mathbf B}-\l I)^{-1}$ is compact
 in $L^2_\rho$.
\end{lemma}

By Lemma \ref{lemspec}, the   centre and stable subspaces of $\BB$
are given by
$
E^c = {\rm Span}\{\psi_0= F\}$ and $E^s = {\rm Span}\{\psi_\b, \,
|\b|>0\}$.

\ssk

Consider next the adjoint (in the dual $L^2$-metric) operator
\begin{equation}
\label{B2} {\bf B}^* = -(-\D)^m -  \textstyle{\frac 1{2m}}\, y
\cdot \nabla \, .
\end{equation}
 For $m \ge 2$, we  treat  $\BB^*$ in
$L^2_{\rho^*} $ with the exponentially decaying weight function
 \begin{equation}
\label{rho2}
 \tex{
 \rho^*(y) = \frac 1{\rho(y)} \equiv {\mathrm
e}^{-a|y|^\a} > 0.
 }
\end{equation}

\begin{lemma}
\label{lemSpec2} {\rm (i)}  $\BB^*: H^{2m}_{\rho^*}  \to
L^2_{\rho^*} $ is a bounded linear operator with the same spectrum
$(\ref{spec1})$ as $\BB$.
Eigenfunctions $\psi_\b^*(y)$ with $|\b|=l$ are $l$th-order
generalized Hermite  polynomials given by
\begin{equation}
\label{psidec}
 \psi_\b^*(y) =
 \textstyle{
  \frac {1}{ \sqrt{\b !}} \, \Big[ y^\b +
\sum\limits_{j=1}^{[ |\b|/2m ]} \frac 1{j !}(-\Delta)^{m j} y^\b
\Big].
 }
\end{equation}

{\rm (ii)} The set $\Phi^*=\{\psi_\beta^*\}_{|\b| \ge 0}$ is
complete and  resolvent $({\mathbf B}^*-\l I)^{-1}$   compact in
$L^2_{\rho^*}$.

\end{lemma}

It follows that the orthonormality condition holds
 \be
 \label{Ortog}
\langle \psi_\b, \psi_\g^* \rangle = \d_{\b\g},
 \ee
where $\langle \cdot, \cdot \rangle$ denotes the standard
$L^2(\ren)$ inner product and $\d_{\b \g}$ is Kronecker's delta.

Using (\ref{Ortog}), we introduce the subspaces of eigenfunction
expansions and begin with the operator $\BB$. We denote
 by $\tilde
L^2_\rho$ the subspace of eigenfunction expansions $v= \sum c_\b
\psi_\b$ with coefficients $c_\b = \langle v, \psi^* \rangle$
defined as the closure of the finite sums $\{\sum_{|\b| \le M}
c_\b \psi_\b\}$ in the norm of $L^2_\rho$. Similarly, for the
adjoint operator $\BB^*$, we define the subspace $\tilde
L^2_{\rho^*} \subseteq L^2_{\rho^*}$.
 Note that since the
operators are not self-adjoint and the eigenfunction subsets are
not orthonormal, in general, these subspaces can be different from
$
 L^2_{\rho}$ and $L^2_{\rho^*}$, and
the equality is guaranteed in the self-adjoint case $m=1$,
$a=\frac 1 4$ only.


\subsection{\bf Existence of similarity profiles close to
transcritical  bifurcations}
 \label{Sect6}

Consider the elliptic problem (\ref{1.7R}).
Using the above Hermitian spectral analysis of the operator pair
$\{\BB,\BB^*\}$, we formulate the bifurcation problems, which
guarantee the existence of a similarity solution in a
neighbourhood of bifurcation points. In fact, our consideration is
quite similar to that for the second-order case in \cite{GMPOscI},
so we may omit some details. Since $p<p_0$, our analysis is
 performed in the
subcritical Sobolev range:
 \be
   \label{Sob1}
    \mbox{$
   1<p<p_S= \frac {N+2m}{N-2m} \LongA H^{m}_\rho(\ren) \subset
   L^{p+1}_\rho(\ren) \,\,\, \mbox{compactly}.
    $}
    \ee


Taking $p$ close to the critical values, as defined in (\ref{h2}),
we look for small solutions of (\ref{1.7R}).  At $p = p_l$, the
linear operator $\BB_1$ has a nontrivial kernel, hence:


\begin{proposition}
\label{PrBif}
 Let for an $l \ge 0$, the
eigenvalue $\l_l = -\frac l{2m}$ of operator $(\ref{ODEf})$ is of
 odd multiplicity. Then the critical exponent $(\ref{h2})$
   is a bifurcation point for
the  problem $(\ref{1.7R})$.
\end{proposition}

\noindent {\em Proof.}
 Consider in $L^2_\rho $ our equation written as
  \be
  \label{treq1m}
 \hat{\BB} f = -(1+c_1)f - |f|^p,\quad \mbox{where} \,\,\, \hat{\BB} =
\BB_1 - (1+c_1)I \equiv \BB-I.
 \ee
It follows that the spectrum $\s(\hat{\BB})=\{-1-
 \frac l{2m}\}$ consists of strictly negative eigenvalues. The
inverse operator $\hat{\BB}^{-1}$ is known to be compact,
\cite[Prop.~2.4]{Eg4}. Therefore, in
 the corresponding integral equation
 \be
 \label{Inteqm}
 f = \hat {\bf A}(f) \equiv - (1+c_1) \hat {\BB}^{-1} f -
\hat{\BB}^{-1} (|f|^p),
 \ee
the right-hand side is a compact Hammerstein operator; see
\cite[Ch.~V]{Kras} and applications in \cite{BGW1, GHUni, GW2}. In
view of the known spectral properties of $\hat {\bf B}^{-1}$,
bifurcations   in the problem (\ref{Inteqm}) occur if the
derivative $\hat {\bf A}'(0)=- (1+c_1) \hat{\BB}^{-1}$ has the
eigenvalue $1$ of odd multiplicity, \cite{KrasZ, Kras}. Since
$\s(\hat{\bf A}'(0))= \{(1+c_1)/(1+\frac l{2m})\}$, we obtain the
critical values (\ref{h2}). By construction, the solutions of
(\ref{Inteqm}) for $p \approx p_l$ are small in $L^2_\rho $ and,
as can be seen from the properties of the inverse operator, $f$ is
small in the domain   $H^{2m}_\rho $ of ${\bf B}$.
 Since the weight (\ref{rho44}) is a monotone growing function as $|y| \to
\infty$, using the known asymptotic properties of solutions of
(\ref{1.7R}), $f \in H^{2m}_\rho $ is a uniformly bounded,
continuous function (for $N<2m$, this directly follows from
Sobolev's embedding theorem).
$\qed$

\ssk

 Thus, $l=0$ is always a bifurcation point since $\l_0=0$ is simple.
 In general, for $l=1,2,...$ the odd multiplicity occurs depending on the
 dimension $N$.
For instance, for $l=1$, the multiplicity is $N$, and, for $l=2$,
it is $\frac {N(N+1)}2$.
  In the case of even multiplicity of
 $\l_l$, an extra analysis is necessary to guarantee that a
 bifurcation occurs \cite{KrasZ} using the rotation $\g_l$ of the vector
 field
  corresponding to the nonlinear term in (\ref{Inteqm})
  on the unit sphere in the eigenspace $\Phi_l= {\rm Span}\{\psi_\b, \,
 |\b|=l\}$ (if $\g_l \not = 1$, then bifurcation occurs).
  We do not perform this study here and note that the
  non-degeneracy of this vector field is not straightforward; see
  related comments below.
It is crucial that, for
 main applications, for $N=1$ and for the radial setting in $\ren$,
 the eigenvalues
(\ref{spec1}) are simple and (\ref{h2}) are always  bifurcation
points.
Unlike Proposition 3.2
 for $m=1$ \cite[\S~3]{GMPOscI}, we have
the following result describing the local behaviour of bifurcation
branches occurring in the main applications, see \cite{Kras} and
\cite[Ch.~8]{KrasZ}. Unlike the case $m=1$, some bifurcations
become {\em transcritical}.

\begin{proposition}
\label{PBif2m}
 Let $\l_l$ be a simple eigenvalue of $\BB$ with  eigenfunction
 $\psi_l$, and let
 \be
  \label{kap1m}
\kappa_l = \langle |\psi_l|^p, \psi_l^* \rangle \not = 0.
 \ee
 Then the $p$-bifurcation branch
 crosses transversely the $p$-axis at $p=p_l$.

\end{proposition}

We next describe the behaviour of solutions for $p \approx p_l$
and apply the classical Lyapunov--Schmidt method, \cite[Ch.~8]
{KrasZ}, to equation (\ref{Inteqm}) with the operator $\hat {\bf
A}$ that is differentiable at $0$. Since, under the assumptions of
Proposition \ref{PBif2m}, the kernel $E_0={\rm ker\,}\hat {\bf
A}'(0)={\rm Span\,}\{\psi_l\}$ is one-dimensional, denoting by
$E_1$ the complementary (orthogonal to $\psi_l^*$) invariant
subspace, we set
  \be
  \label{EXP1}
  \tex{
  f= F_0+F_1 \whereA
 F_0= \e_l \psi_l \in E_0 \andA
F_1= \sum_{k \not = l} \e_k \psi_k \in E_1.
 }
 \ee
  Let $P_0$ and $P_1$,
$P_0+P_1=I$, be projections onto $E_0$ and $E_1$ respectively.
 Projecting (\ref{Inteqm})
  onto $E_0$ yields
 \be
 \label{eqIntm}
\g_l \e_l = - \langle \hat \BB^{-1} (|f|^p), \psi_l^* \rangle,
 \quad \g_l =1- \textstyle{\frac{1+c_1}{1+l/2m}}= \textstyle{\frac{(N+l)s}{(p-1)(2m+l)}},
 \quad s=p-p_l.
 \ee
By bifurcation theory (see \cite[p.~355]{KrasZ} or
\cite[p.~383]{Deim}, where $\hat {\bf A}'(0)$ is Fredholm of index
zero), $F_1 = o(\e_l)$ as $\e_l \to 0$, so that $\e_l$ is
calculated from (\ref{eqIntm}) as:
 \be
 \label{Compmm}
  \tex{
\g_l \e_l = -|\e_l|^p \langle \hat \BB^{-1} |\psi_l|^p, \psi_l^*
\rangle
 + o(\e_l^p) \,\, \Longrightarrow \,\,
 |\e_l|^{p-2}\e_l = \hat c_l (p-p_l)[1+
 o(1)],
 }
  \ee
   where $\hat c_l = \textstyle{\frac{
(N+l)^2}{4m^2 \kappa_l}}$.
   We have used the following
calculation:
  $$
   \tex{
  \langle \hat \BB^{-1} |\psi_l|^p, \psi_l^* \rangle =
 \langle |\psi_l|^p,  (\hat \BB^*)^{-1} \psi_l^* \rangle =
 -\frac{\kappa_l}{1+ \frac l{2m}}.
 }
 $$
Recall the identity $(\hat \BB^{-1})^*=(\hat \BB^*)^{-1}$.

It follows from the algebraic equation in \ef{Compmm} that the
bifurcations are transcritical provided that $\kappa_l \not =0$,
while the sign of $\kappa_l$ determines how the branches cross the
$p$-axis.
 For $l=0$, it is easy to see that $\kappa_0>0$, since $\psi_0=F$ and
$\psi_0^*=1$, so that
 \be
 \label{kap00}
 \mbox{$
  \kappa_0= \langle |\psi_0|^p, \psi_0^* \rangle = \int |F|^p>0 \quad (p=p_0).
  $}
  \ee
Moreover, for $p_0 \approx 1^+$,
 by the definition of $F$ in \ef{ODEf}, we have that $\kappa_0
 \approx 1$.
 The positivity or negativity of the scalar product
(\ref{kap1m}) for $l \ge 1$ and  arbitrary $p>1$ is not
straightforward, and we should rely on a delicate numerical
evidence; see \cite{GHUni}. It turns out that $\kappa_l$ can be
both positive and negative for different $l \ge 1$.

Let us note the following principal difference in comparison with
the case of monotone nonlinearity $+|f|^{p-1}f$ studied before for
$m=1$, \cite{GMPOscI}. It turns out that
\be
\label{lodd}
 \kappa_l=0 \quad \mbox{for odd $l=1,3,5,...$},
 \ee
 since in
\ef{kap1m}, $|\psi_l(y)|^p$ is even, while the polynomial
$\psi_l^*(y)$ is odd. This means the following:
 \be
 \label{Bifl1}
 p=p_l \quad \mbox{for odd $l \ge 1$ are not {\bf ``standard"} bifurcation points.}
  \ee
  The corresponding ``non-standard" bifurcation phenomenon will be
  discussed shortly.
On the other hand, for the standard monotone nonlinearity as in
\ef{2mGH}, all critical exponents $\{p_l\}$ are pitchfork
bifurcations, \cite{GHUni}.

\ssk

Thus, under the  assumption $\kappa_l \not =0$
on the coefficients (\ref{kap1m}), we obtain a countable sequence
of  bifurcation points (\ref{h2}) satisfying $p_l \to 1^+$ as $l
\to \infty$, with typical transcritical  bifurcation branches
appearing in a neighbourhood. The behaviour of solutions in
$H^{2m}_\rho$ and uniformly in $\ren$,  for $p<p_S$, takes the
form
 \be
  \label{Vexpnm}
 f_l(y) = |\hat c_l (p-p_l)|^{\frac 1{p-1}} \, {\rm
 sign}\,(\hat c_l (p-p_l))\,
( \psi_l(y) + o(1)) \quad \mbox{as} \,\,\, p \to p_l.
 \ee
Instability of all these local branches of similarity profiles is
studied similar to the case $m=1$ in \cite[\S~2]{GMPOscI};
 see also \cite{GHUni}.

\subsection{Lyapunov--Schmidt branching equation in the general multiple case: non-radial patterns}

Let now $\l_l= - \frac l{2m}$ have multiplicity $M=M(l) > 1$ given
by the binomial coefficient
 \be
   \label{93}
    \tex{
  M(l)=  {\rm dim} \,W^c(\BB - \l_l I)=C_{N+l-1}^l= \frac
    {(N+l-1)!}{l!(N-1)!}, \quad \mbox{so that}
     }
     \ee
 \be
 \label{e1}
E_0= {\rm ker}(\BB- \l_l I)= {\rm
Span}\{\psi_{l1},...,\psi_{lM}\}.
  \ee
Then,
 looking for a solution
 \be
 \label{e2}
 f=f_0+f_1 \withA f_0= \e_1 \psi_{l1}+...+\e_M \psi_{lM} \whereA f_1
 \bot E_0,
  \ee
  and substituting into the equation (\ref{Inteqm}), multiplying by $\psi_{li}^*$, and
  denoting,
  as usual, $s=p-p_l
 $, $0<|s| \ll 1$, we obtain the following {\em generating} system of $M$ algebraic equations:
\be
 \label{e3m}
  \mbox{$
   \e_i= \frac {2m} {s(N+l)^2}\, \int  |\e_1 \psi_{l1}+...+\e_M
 \psi_{lM}|^{p}\,  \psi_{li}^* \equiv
 D_i(\e_1,...,\e_M), \quad i=1,2,...,M.
  $}
  \ee
   Here $p=p_l$.
Denoting $x=(\e_1,...,\e_M)^T \in \re^M$, the system (\ref{e3m})
is written as a fixed point problem for the given nonlinear
operator
 ${\bf  D}=(D_1(x),...,D_M(x))^T$,
 \be
 \label{e91}
  x= {\bf D}(x) \inB \re^M.
   \ee
In the second-order case $m=1$ \cite[\S~3]{GMPOscI}, the system
\ef{e91} was variational, that allowed us to get a multiplicity
result.
 In view of the dual metric in \ef{Ortog}, for
any $m \ge 2$, the algebraic system \ef{e3m} is not variational,
so the multiplicity of admissible solutions remains an open
problem.

\ssk

Global extension of the above local $p$-bifurcation branches is
performed by classic theory of nonlinear compact operators,
\cite{Deim, Kras, KrasZ}. However, since the problem is not
variational, nothing prevents existence of closed $p$ sub-branches
or appearances of turning, saddle-node bifurcations (we will show
that this can actually happens), so that the total number and
structure of solutions for any $p \in (1,p_0]$
 remain a difficult problem. We
 will then inevitably should rely on careful numerics.


 \subsection{``Non-standard" pitchfork bifurcations for $\kappa_l=0$}
  \label{S2.5}

Without loss of generality, we consider the simplest case $l=1$,
$N=1$, $m=2$ (then $p_1=1+m=3$ by \ef{h2}), where, from \ef{kap1m}
and \ef{eigen}, \ef{psidec}, it is clear that $\kappa_1$ vanishes:
 \be
 \label{kap1}
 \tex{
  \kappa_1=\langle |\psi_1|^{p_1},\psi_1^* \rangle \equiv \int_\re
  \big |F'(y)|^{3} y \, {\mathrm d}y=0.
  }
  \ee
 Next, unlike the standard approximation \ef{EXP1} close to $p=3$,
 we now use an improved one given by the expansion on the 2D
 invariant subspace $E_{12}={\rm Span}\{\psi_1,\psi_2\}$ (this choice
 will be explained below):
  \be
  \label{b12}
  f= F_{12} + F_3 \whereA F_{12}=\e_1 \psi_1+\e_2 \psi_2 \andA F_3
  \bot E_{12},
   \ee
   with the scalar parameters $\e_1$, $\e_2$ to be determined.
   For simplicity, we next use the differential version of the
   integral equation \ef{Inteqm} for
   $l=1$:
    \be
    \label{dif1}
     \tex{
 (\BB - \l_1 I)f - \frac s4 \, f=|f|^p+... \whereA s=p-p_1
 }
 \ee
 and where we omit the $O(s^2)$-term.
   Substituting \ef{b12} into \ef{dif1} and
   projecting onto corresponding one-dimensional subspaces, quite
   similar to the system \ef{e3m}, we obtain the following
   asymptotic system of two algebraic equations:
    \be
    \label{sys12}
     \left\{
      \begin{matrix}
     - \frac s4 \, \e_1= -\int |\e_1 \psi_1+\e_2 \psi_2|^3 \psi_1^*
     +...,\quad \,\,\, \ssk\ssk\\
     - \frac 14 \, \e_2 =
      \frac s4 \, \e_2- \int |\e_1 \psi_1+\e_2
     \psi_2|^3
     \psi_2^*+...\,,
      \end{matrix}
      \right.
      \ee
   where we put $\l_2=- \frac l{2m}=-\frac 12$ (for $l=m=2$) and where we have omitted higher-order terms associated with the
   orthogonal $F_3$ in \ef{b12} and via replacing $p$ by $p_1=3$ in the integrals on the right-hand sides.
Then, the second equation, as $s \to 0$, gives the dependence of
$\e_2$ on the leading expansion coefficient $\e_1$ on $E_{12}$:
 \be
 \label{eee2}
  \tex{
  \e_2=2 |\e_1|^3 \mu_{12}+... \whereA \mu_{12}=\int
  |\psi_1|^3\psi_2^* \ne 0 \andA \psi_2^*= \frac 1{\sqrt 2}\,
  y^2.
  }
  \ee
  It is crucial that, unlike in \ef{kap1}, the coefficient
  $\mu_{12}$ is given by the integral of some {\em even}
  function, so that, now, the assumption $\mu_{12} \ne 0$ is not
  that restrictive and can be quite reliably checked numerically.

Next, the first equation in \ef{sys12}, after a Taylor expansion
in the integral, by using that $\e_2 = o(\e_1)$ as $s \to 0$,
provides us with the necessary bifurcation scalar equation on
$\e_1$,
 \be
 \label{eee11}
 \tex{
\frac s4 \, \e_1= \int|\e_1 \psi_1|^3 \psi_1^* +3 \e_1^2 \,{\rm
sign}\, \e_1 \e_2 \nu_{12}+..., \,\, \mbox{with}\,\, \nu_{12}=
\int \psi_1^2 ({\rm sign}\, \psi_1)\, \psi_2 \psi_1^* \ne 0, }
  \ee
   where, again, in
$\nu_{12}$, we face an even function in the integral. Since the
first coefficient vanishes by \ef{lodd} for $l=1$,  $\int
|\psi_1|^3\psi_1^*=0$, using the dependence \ef{eee2}, we obtain
 \be
 \label{ee11}
  \tex{
  |\e_1|^5= \hat c_{12} \, s +... \whereA \hat c_{12}= \frac 1{24
  \mu_{12}\nu_{12}}.
  }
   \ee
 It follows
  that
 we thus deal with a {\em pitchfork} bifurcation at $p=p_1=3$, which is
 subcritical if $\hat c_{12}<0$ and supercritical if $\hat c_{12} >0$.

Overall, the bifurcation branches take the following form: for,
e.g., $\hat c_{12}>0$,
 \be
 \label{bif12}
  \tex{
  f(y)= \pm [\hat c_{12}(p-3)]^{\frac 15} \psi_1(y)+ 2\mu_{12} [\hat c_{12}(p-3)]^{\frac 35} \psi_2(y)
 +... \asA p \to 3^+.
 }
 \ee
 We will reveal this kind of bifurcation numerically in Section
 \ref{SectN4}.
 Note that this non-standard bifurcation branch near $p=3$
 is more ``steep", $\sim O((p-3)^{\frac 15}))$, than the standard one in
\ef{Vexpnm}, which, for $l=1$, is of the order $\sim
O(\sqrt{p-3})$.


\ssk

One can see that a similar bifurcation scenario, under the
vanishing assumption \ef{lodd}, can be developed by using other
invariant subspaces rather than that in \ef{b12}. The crucial
conditions then remain the same: the corresponding coefficients
$\mu_{..}$ in \ef{eee2} and $\nu_{..}$ in \ef{eee11} must be
non-zero, which is possible by mixing even and odd eigenfunctions
in the subspace, depending on the multiindices chosen. This has an
interesting and surprising consequence:
 \be
 \label{cons1}
 \mbox{(\ref{lodd}): there can be more than one
 bifurcation branch, even for 1D eigenspace.}
  \ee
In Section \ref{SectN4}, we will observe this numerically for the
one-dimensional eigenspace.

\subsection{Transversality of intersections of subspaces}
 \label{S.6.3}

This was a permanent subject of an intense study for nonlinear
second-order parabolic equations; see related key references and
further  comments in \cite[\S~6.2]{GMPOscI}. We briefly recall
these important results below.
 Namely, this problem was completely  solved  rather recently for a scalar
 reaction-diffusion equation on a circle of the form
 \be
 \label{RD11}
 u_t={\bf A}(u) \equiv u_{xx}+ g(x,u,u_x), \quad x \in S^1= \re/2\pi{\mathbb
 Z},
  \ee
  where the nonlinearity $g(\cdot)$ satisfies necessary conditions for existence of
  global classical bounded solutions for arbitrary bounded smooth
  initial data.
 Then, if $f$ is a   hyperbolic {\em equilibrium} of ${\bf A}$, ${\bf
A}(f)=0$, known to be generic (or a {\em rotating wave}), then the
global stable and unstable subspaces of ${\bf A}'(f)$ span the
whole functional space
 $X^\a=H^{2 \a}(S^1)$, $ \a \in \big( \frac 32, 1)$, where the
 global semiflow is naturally defined, i.e.,
  \be
  \label{RD12}
   W^{s}({\bf A}'(f))  \oplus  W^{u}({\bf A}'(f))= X^\a,
    \ee
    so that these subspaces {\em intersect transversely}.
 It is crucial that such a complete analysis can be performed in
 1D only, since it is based on Sturmian zero set arguments (see
 \cite{GalGeom} for main references and various extensions of
 these fundamental ideas), so, in principle, cannot be extended
 to equation in $\ren$.
    We refer to most recent papers \cite{Czaja08, FiedII08,  Joly10},
    where earlier key references and most advances results on the
    transversality and connecting orbits can be found.



We perform our transversality analysis for $p$ close to the
bifurcation points $ p \approx p_l$ in \ef{h2} by using
bifurcation theory from Section \ref{Sect6}:

\begin{proposition}
 \label{Pr.Trans}
  Fix, for a given $p \approx p_l$, $p \ne p_l$,
 a hyperbolic equilibrium $f_\b$, with a $|\b|=l$, of the operator ${\bf A}$ in
 $\ef{1.7R}$,
 \be
 \label{C27m}
  \tex{
  {\bf A}(f)=-(-\D)^m f + \frac 1{2m} \, y \cdot \n f+ \frac 1{p-1}
  \, f + |f|^p.
   }
   \ee
  Then the transversality conclusion holds:
 \be
  \label{RD13}
   W^{s}({\bf A}'(f_\b))  \oplus  W^{u}({\bf A}'(f_\b))= H^{2m}_\rho(\ren).
    \ee
\end{proposition}

\noi{\em Proof.} 
 It follows from \ef{C27m}
 and the expansion
\ef{Vexpnm} that, for $p= p_l + \e$, with $0<|\e| \ll 1$,
 \be
 \label{91}
  \begin{matrix}
  {\bf A}'(f_\b)=-(-\D)^m  + \frac 1{2m} \, y \cdot \n+  \frac 1{p-1}\, I
  + p|f_\b|^{p-1} \, {\rm sign}\, f_\b\, I
 \ssk\ssk \\
   = (\BB - \l_l I)+ p|\hat c_l|\, |\e| |\psi_l|^{p-1} {\rm
   sign}\,(c_l \e\psi_\b)+... \quad (l=|\b|) \, .
   \end{matrix}
   \ee
Therefore, for $p=p_l$, the following analogy of \ef{RD13} is
valid:
 \be
 \label{92}
 {\bf A}'(f_\b)=\BB - \l_l I \LongA W^s(\BB - \l_l I) \oplus W^u(\BB -
 \l_l
 I) \oplus W^c(\BB - \l_l I)= H^{2m}_\rho(\ren),
  \ee
  and  ${\rm dim} \,W^c(\BB - \l_l I)$ is
  equal to the algebraic multiplicity \ef{93} of $\l_l= - \frac
  l{2m}$.
 By the assumption of the hyperbolicity of $f_\b$ and in view of
 small perturbations (see, e.g., \cite{BS, Kato})
  of all the eigenfunctions of ${\bf A}'(f_\b)$ for any $|\e| \ll 1$, $\e
 \not =0$, which remain complete and closed as for $p=p_l$, we
 arrive at \ef{RD13}. Recall that, since by \ef{91},  ${\bf
 A}'(f_\b)$, with eigenfunction $\{\hat\psi_\b\}$, is a small perturbation of $\BB-\l_\b I$
 (with eigenfunctions $\{\psi_\b\}$) and, in addition, the
 perturbation is exponentially small as $y \to \iy$, the
 ``perturbed" eigenfunctions $\hat \psi_\b(y)$ remain a small perturbation of the
 known $\psi_\b(y)$ in any bounded ball, and sharply approximate
 those as $y \to \iy$. Therefore, close to $p=p_l$, there is no
 doubt that the well-known condition of completeness/closure of $\{\hat\psi_\b\}$ (the
 so-called {\em property of stability of the basis})  is, indeed, valid:
  $$
   \tex{
    \sum_{(\b)} \|\psi_\b\|_\rho \, \| \hat \psi_\b-\psi_\b \|_\rho
    <1. \quad \qed
    }
    $$


\ssk

  Thus, close to any bifurcation point $p=p_l$, we precisely know both
   the dimensions of the unstable subspace of ${\bf
  A}'(f_\b)$ of any hyperbolic equilibrium $f_\b$ (and, sometimes,
  we can prove the latter) and  the
  corresponding eigenfunctions $\{ \hat \psi_\b\}$:
 \be
 \label{96}
 \tex{
 \mbox{ by continuity, for all} \,\,\,p \approx p_l: \quad \hat \l_l \approx -\l_l=
 \frac l{2m} \andA \hat \psi_\b \approx \psi_\b,
 }
 \ee
 where convergence of eigenfunctions as $p \to p_l^-$ holds in $H^{2m}_\rho$
 and uniformly in $\ren$.

 Furthermore, moving along the given bifurcation $p$-branch, the
 transversality persists until a saddle-node bifurcation appears, when a centre
 subspace
 for ${\bf A}'(f_\b)$ occurs, and hence \ef{RD13} does not apply.
 If such a ``turning" point of a given $p$-branch  does not appear (but sometimes it does;
 see Section \ref{SectN4} below),
  the transversality
persists globally in $p$.

\section{Numerical results: extension of even $p$-branches}
 \label{S3}



Thus, the above bifurcation analysis establishes existence of a
countable set of transcritical $p$-bifurcations at $p=p_l$ for
even $l$.
 As we have mentioned, since \ef{1.7R} is
not variational for  $m \ge 2$, we do not have any chance to use
power tools of category-genus-fibering theory in order to
guarantee nonlocal extensions of $p$-branches of similarity
profiles $f(y)$. However, as is well-known from compact nonlinear
integral operator theory \cite{Deim, Kras, KrasZ}, these branches
are always extensible, but can end up at other bifurcation points,
so their global extension for all $p>p_l$ is not straightforward.
Actually, we show that precisely this happens for $m=2$ in 1D.

\subsection{Preliminaries for $m=2$: well-posed shooting of even profiles}

 We first concentrate on the simplest fourth-order case:
 \be
 \label{F1}
  \tex{
 N=1 \andA m=2, \quad \mbox{so that} \quad p_0=1 + \frac {2m} N=5,
 }
  \ee
in order to exhibit typical difficult and surprising behaviours of
global $p$-branches of the first similarity profile $f_0(y)$,
which bifurcates from the first critical exponent $p_0=5$ in
\ef{F1}.
 We also compare $f_0$ in dimensions $N=1$, 2, 3, and 4.
 For convenience, we will denote by $f_l(y)$ the profiles
that bifurcate at the corresponding critical $p_l$ and hence, by
\ef{Vexpnm},  ``inherit" the nodal set structure of the
eigenfunction $\psi_l(y)$ in \ef{eigen} for $N=1$.

In the case \ef{F1}, the problem \ef{1.7R} becomes an ODE one:
 \be
  \label{09Nm}
   \left\{
   \begin{matrix}
 {\bf A}(f) \equiv -f^{(4)} + \frac 14 \,y f'  + \frac 1{p-1}\, f
 +
 |f|^{p}=0 \forA y>0,\ssk\ssk \\
 f(y) \,\,\, \mbox{decays exponentially fast as $y \to \pm \iy$}.
 \qquad\qquad\,\,\,
  \end{matrix}
   \right.
    \ee

We first easily prove the following result, somehow confirming our
bifurcation analysis:

\begin{proposition}
 \label{Pr.Non}
{\rm (i)} In the critical case $p=5$, the only solution of
$\ef{09Nm}$ is
 $f=0$; and

{\rm (ii)} The total mass of solutions of $\ef{09Nm}$ satisfies
 \be
 \label{Int12}
 \tex{
 \int f <0 \forA p<p_0=5 \andA \int f >0 \forA p>5.
 }
 \ee

  \end{proposition}

  \noi{\em Proof.} Integrating the ODE \ef{09Nm} over $\re$ yields the
  following identity:
   \be
   \label{Int11}
    \tex{
     \int  |f|^p=  \frac {p-5}{4(p-1)} \, \int f.
     \qed
     }
     \ee


\noi{\bf Remark on bifurcation analysis.} Firstly, according to
the bi-orthogonality \ef{Ortog},
 \be
 \label{tt1}
 \tex{
 \int \psi _l=0 \quad \mbox{for all} \quad l=1,2,3,...\, , \andA
 \int \psi_0=\int F=1,
 }
 \ee
 so we see that \ef{Vexpnm} somehow ``contradicts" \ef{Int12}.
 However, there is no any controversy here: indeed, \ef{EXP1}
 assumes, in 1D, the following expansion:
  \be
  \label{ff1}
  f= \e_l \psi_l+ \e_0 \psi_0+... \, ,
  \ee
  where we keep the only eigenfunction $\psi_0$ with  the unit
  non-zero mass. Then, the identity \ef{Int11} is perfectly valid
  provided that
   \be
   \label{ff2}
    \tex{
    \e_0= \frac {4(p-1)}{p-5} \, |\e_l|^p \int |\psi_l|^p (1+o(1))=
    o(\e_l),
    }
    \ee
    so that this small correction in \ef{ff1} allows one
to keep the necessary non-zero mass on any even $p$-bifurcation
branch.

\ssk

As usual, for the  even profile  $f_0$ (and for $f_2$, $f_4$,...),
since the ODE \ef{09Nm} is invariant under the symmetry reflection
$y \mapsto -y$,
  two
symmetry conditions at the origin are imposed,
 \be
 \label{Symm1m}
 f'(0)=f'''(0)=0 \quad(\mbox{then} \,\,\,f(-y) \equiv f(y)).
 \ee

Let us first reveal a natural ``geometric" origin of existence of
various solutions of the problem \ef{09Nm}, \ef{Symm1m}. This is
important for the present non-variational problem, where we do not
have other standard techniques of its global analysis. It is easy
to see that the ODE in \ef{09Nm} admits 2D bundle of proper
exponential asymptotics as $y \to +\iy$:
 \be
 \label{Asyy}
  \tex{
  f(y) \sim  {\mathrm e}^{-\frac {a_0}2\, y^{4/3}}\big[ C_1
   \cos \big(\frac{a_0 \sqrt 3}2 \, y^{\frac 43}\big) + C_2
 \sin \big(\frac{a_0 \sqrt 3}2 \, y^{\frac 43}\big)\big], \,\,\, a_0= 3
 \cdot 2^{-\frac 83},
  }
  \ee
  where $C_{1,2} \in \re$ are arbitrary constants.
 Obviously, \ef{09Nm} also admits a lot of solutions with much slower algebraic
 decay,
  \be
  \label{Bl4}
  f(y) \sim C_0 y^{-\frac 4{p-1}} \asA y \to + \iy, \quad C_0 \in
  \re, \quad C_0 \ne 0,
   \ee
  but these should be excluded from the consideration, so we always must take $C_0=0$.

 These {\bf two} parameters $C_{1,2}$ in \ef{Asyy} are  used to satisfy
  (to ``shoot") also {\bf two} conditions at the origin \ef{Symm1m}.
  Overall, this looks like a well-posed (``{\bf 2--2}", i.e., not over- and under-determined)
   geometric shooting
  problem, but indeed extra difficult ``oscillatory" properties of the ODE
  involved are necessary to guarantee a proper mathematical
  conclusion on existence of solutions and their multiplicity (in
  fact, an infinite number of those).  This will be done with the
  help of numerical methods, and, as was mentioned, the final conclusions are
  striking different from those obtained in \cite{GHUni, GMPSob,
  GMPSobII} for monotone nonlinearities.

  Thus, we arrive at a well posed ``$2-2$" shooting problem:
  denoting by $f=f(y;C_1,C_2)$ solutions having the asymptotic
  behaviour \ef{Asyy} (note that such solutions can blow-up at
  finite $y_0 \ge 0$, but we are interested in those with $y_0(C_1,C_2)<0$; see below), by \ef{Symm1m}, an
  algebraic system of two equations with two unknowns occurs:
   \be
   \label{Sys1}
   \left\{
   \begin{matrix}
   f'(0;C_1,C_2)=0,
 \\
   f'''(0;C_1,C_2)=0.
    \end{matrix}
    \right.
    \ee

 \begin{proposition}
  \label{Pr.Anal}
  For any even integers $p=2,4,6,...$, the system \ef{Sys1} admits
  not more than a countable set of solutions.
   \end{proposition}

   \noi{\em Proof.} For such $p$'s, the ODE \ef{09Nm} has an
   analytic nonlinearity, so by classic ODE theory
   \cite[Ch.~I]{CodL}, both functions in \ef{Sys1} are also
   analytic, whence the result. \quad $\qed$

   \ssk

   We expect that a similar result is true for arbitrary $p>1$,
   but a proof of an analytic dependence on parameters\footnote{As
   is well known, dependence on parameters in such ODE problems
   can be much better than the smoothness of coefficients involved. A
   classic example is: for elliptic operators with just measurable
   coefficients, the resolvent is often a meromorphic function
   of the spectral parameter $\l \in {\mathbb C}$.},
   is expected to be very technical.

\subsection{The first symmetric profile $f_0(y)$}

For solving our problem \ef{09Nm}, we use the {\tt bvp4c} solver
of the {\tt MatLab} with the enhanced accuracy and tolerances in
the range
 \be
 \label{tol1}
 10^{-6} - 10^{-12},
  \ee
  and always, with a proper choice of initial approximations (data),
  observed fast convergence and did not need more than 2000--8000 points, so that each computation usually took
  from 15 seconds to a few minutes.

\begin{figure}
\centering
\includegraphics[scale=0.85]{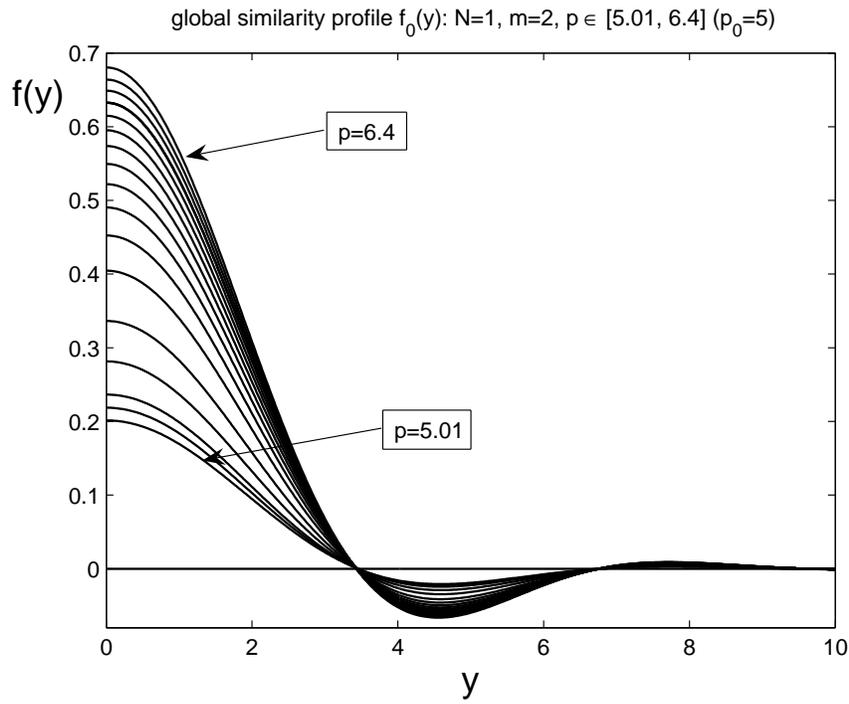}  
\vskip -.5cm \caption{\small The first profile $f_0(y)$ of
(\ref{09Nm}) for  $N=1$, $m=2$ and  $p \in [5.01,6.4]$.}
   \vskip -.1cm
 \label{F00m}
\end{figure}

We begin with Figure \ref{F00m} presenting a general view of the
similarity profile $f_0(y)$ for various $p$ above the critical
exponent $p_0=5$. It is clearly seen that $f_0(y)$ is oscillatory
for large $y$, but definitely has a dominated ``positive hump" on
$y \in (0, 3.4)$, so that overall (cf. \ef{Int12})
 $$
 \tex{
  \int f_0 >0.
   }
   $$
However, by \ef{u01}, this does not imply blow-up of the
corresponding similarity solution $u(x,t)$, since this happens in
the supercritical range $p>p_0=5$, when, in particular, all
sufficiently small solutions are known to be global in time.

\begin{figure}
\centering
\includegraphics[scale=0.8]{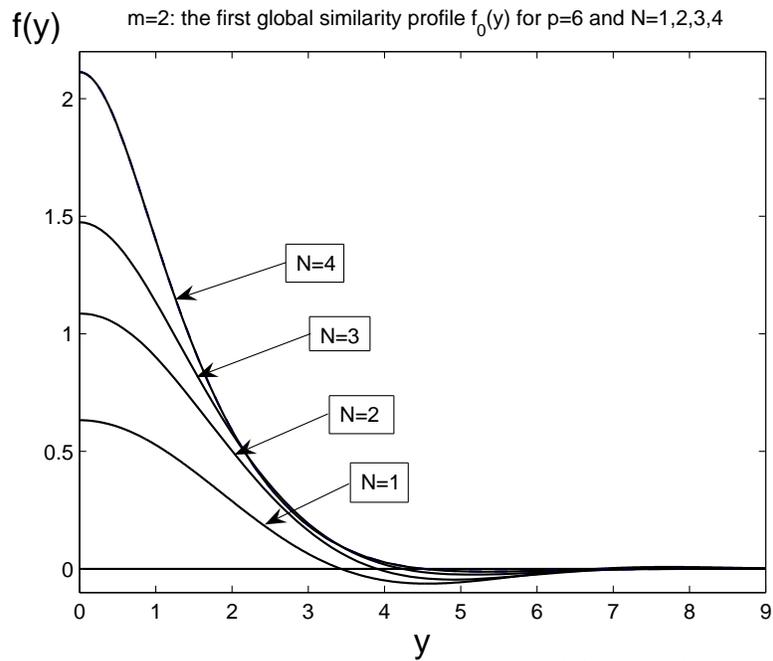}  
\vskip -.5cm \caption{\small The first radially symmetric solution
$f_0(y)$ of (\ref{1.7R}) for $m=2$,  $p=6$, and dimensions $N=1$,
2, 3, 4.}
   \vskip -.1cm
 \label{F00m1N}
\end{figure}

 Figure \ref{F00m1N} shows the
dependence of the radial pattern $f_0=f_0(|y|)$ on dimensions
$N=1$, 2, 3, 4.
 All the profiles look similar and their $L^\iy$-norm, $f_0(0)$,
 increases with $N$. However, the location of the ``positive
 hump" of each $f_0(y)$ remains practically unchanged, as well as the location of the first ``nonlinear
 transversal zero", $y_0 \sim 4$ always; see more below.

\subsection{$p$-branches and further even profiles}

\begin{figure}
\centering
\includegraphics[scale=0.85]{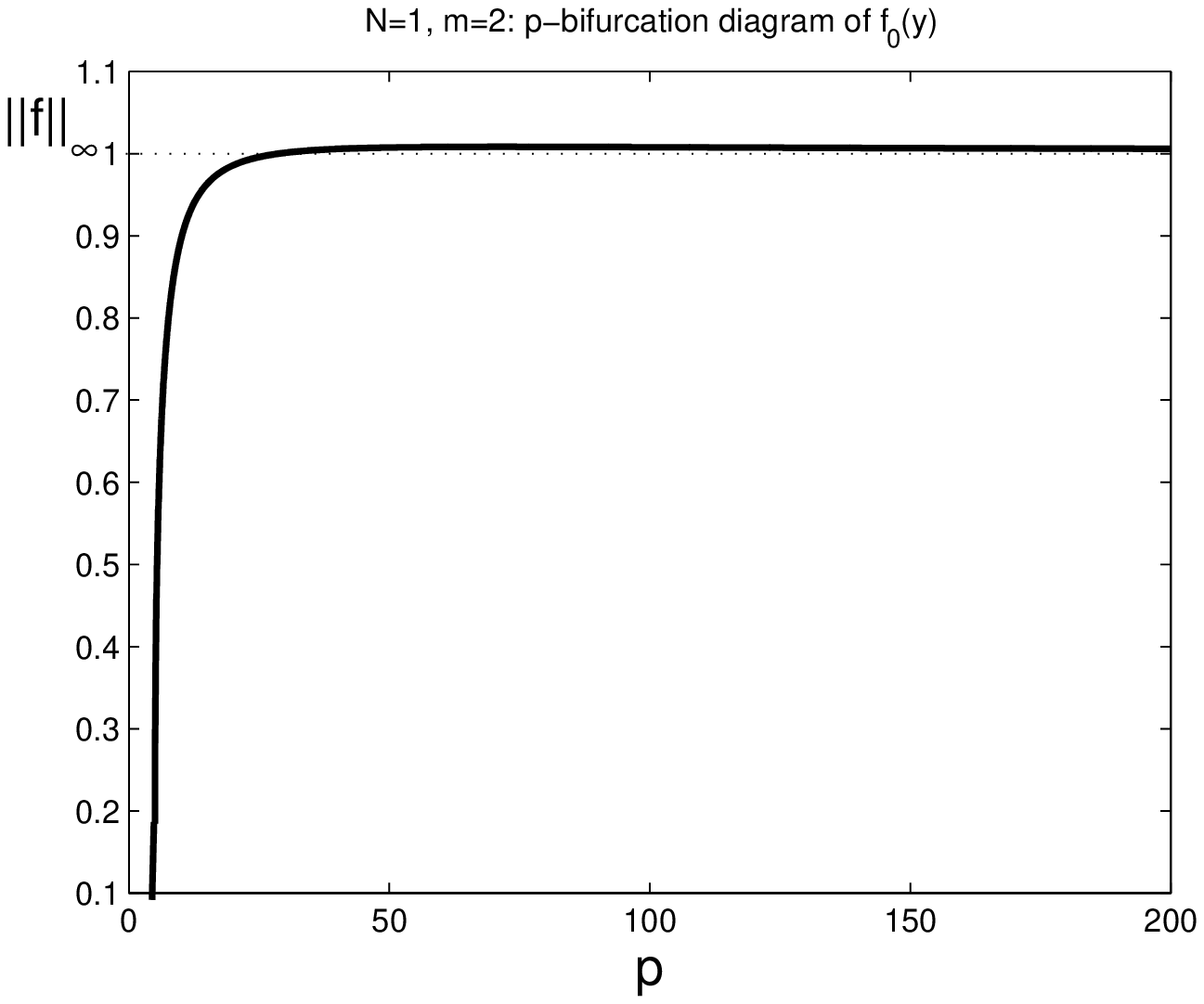}    
\vskip -.5cm \caption{\small $p_0$-branch of $f_0$  for  $N=1$,
$m=2$, extended for $p \in [5.01,200]$.}
   \vskip -.1cm
 \label{F1m}
\end{figure}

\begin{figure}
\centering
\includegraphics[scale=0.85]{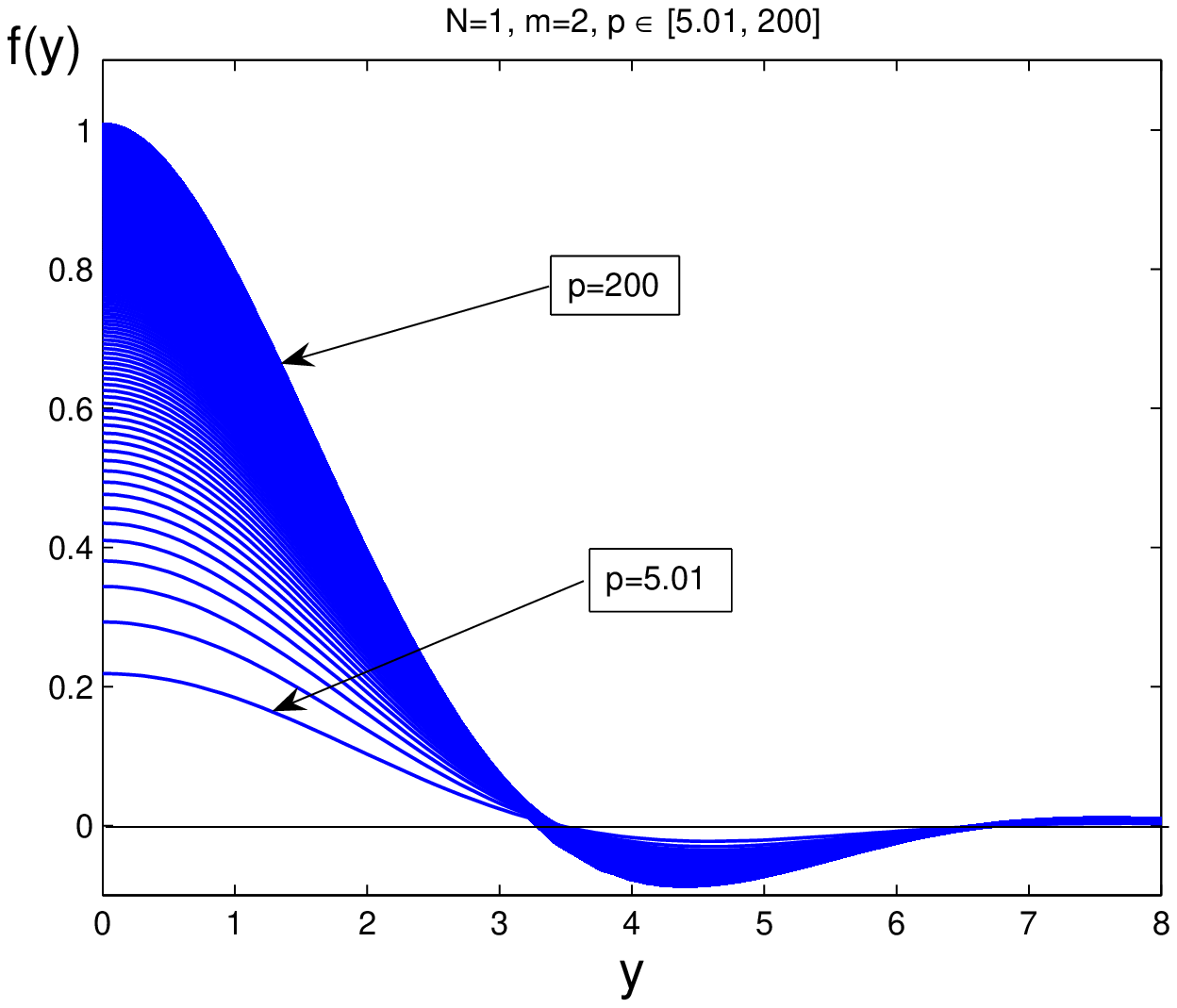}     
\vskip -.5cm \caption{\small $p$-deformation of $f_0$ from Figure
\ref{F1m}; $N=1$, $m=2$, $p \in [5.01,200]$.}
   \vskip -.1cm
 \label{F2m}
\end{figure}

More delicate results are shown
 in Figure \ref{F1m}, where we present the global $p_0$-branch,
 initiated at $p=5^+$ and extended up to $p=200$.
 In particular,  this shows that
  \be
  \label{finf}
  \|f_0\|_\iy \equiv f_0(0) \to 1^+ \asA p \to +\iy,
  \ee
 an asymptotic phenomenon with a possible difficult logarithmically
 perturbed behaviour that was discovered and studied in
 \cite[\S~5]{GHUni} for another model \ef{09Nm} with the monotone nonlinearity $|f|^{p-1}f$.
  The
  deformation of the profile $f_0$ on the same interval $p \in [5.01,200]$  is shown in
  Figure \ref{F2m}, again confirming \ef{finf}.

\begin{figure}
\centering
\includegraphics[scale=0.7]{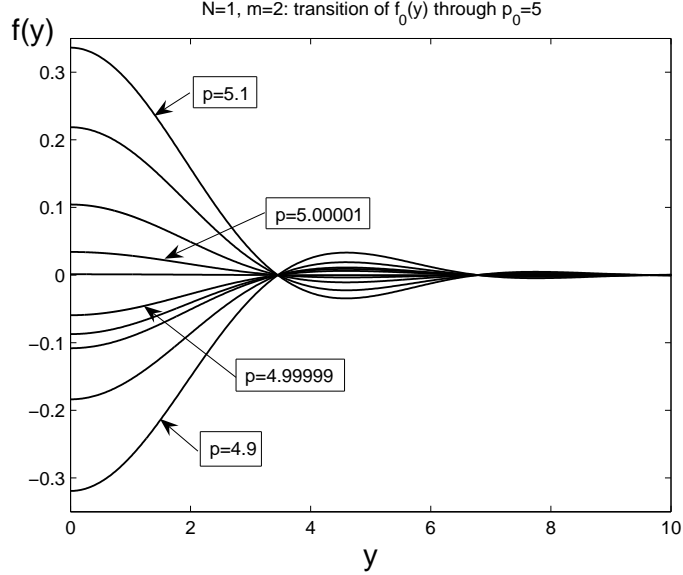}  
\vskip -.5cm \caption{\small Transition of  $f_0(y)$ to $f_2(y)$
of (\ref{09Nm}) for $N=1$, $m=2$ for  $p \approx 5^\pm$.}
   \vskip -.1cm
 \label{F11m}
\end{figure}

We next study the extension of the $p_0$-branch for $p<p_0=5$. The
transition through the first transcritical bifurcation at
$p=p_0=5$ is explained in Figure \ref{F11m}, which shows a clear
spatial similarity of $f_0(y) \sim \pm \psi_0(y)=\pm F(y)$ along
both limits $p \to p_0^\pm=5^\pm$, according to \ef{Vexpnm} for
$l=0$.

\begin{figure}
\centering
\includegraphics[scale=0.75]{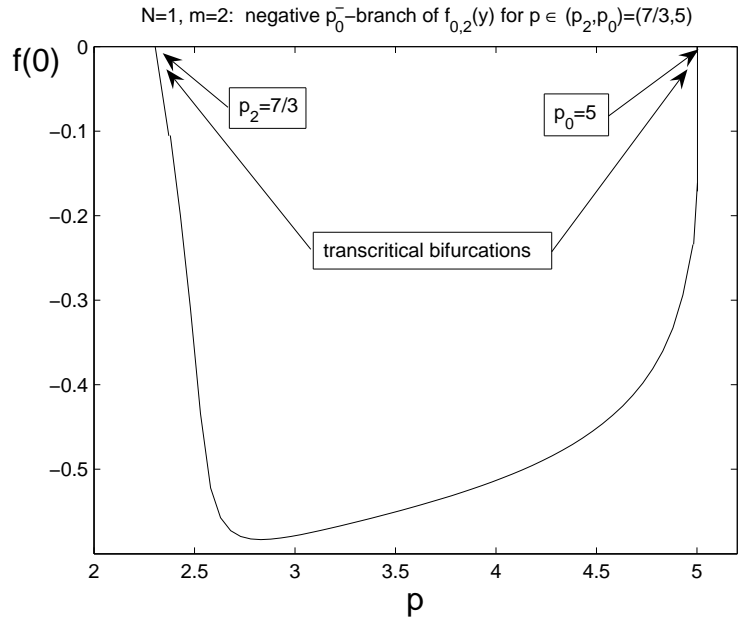}  
\vskip -.5cm \caption{\small $p_2$-branch of $f_2$  for  $N=1$,
$m=2$, extended for $p \in (p_2,p_0)$.}
   \vskip -.1cm
 \label{F3m}
\end{figure}

\begin{figure}
\centering
\includegraphics[scale=0.85]{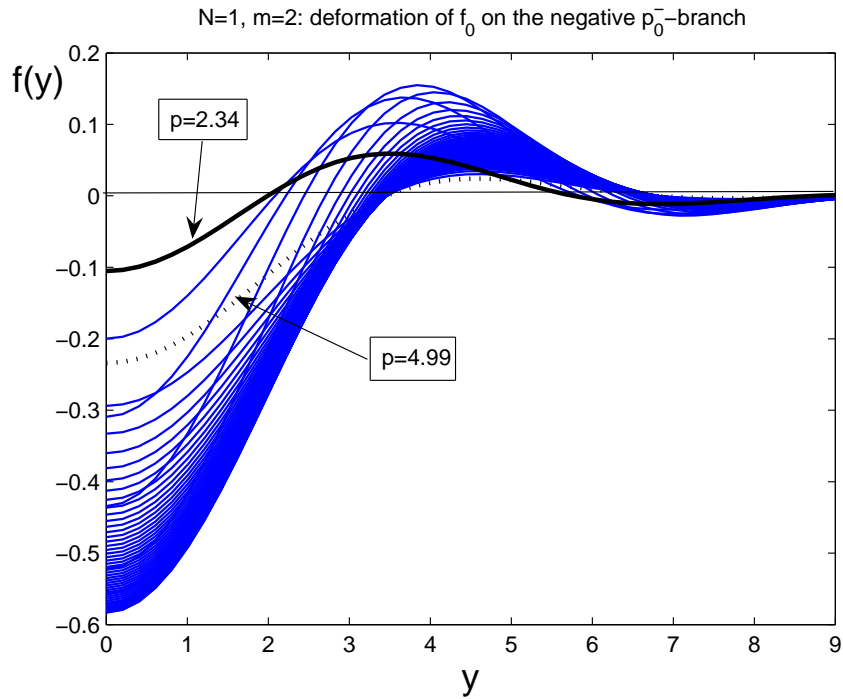}  
\vskip -.5cm \caption{\small $p$-deformation of $f_2$ from Figure
\ref{F3m}; $N=1$, $m=2$, $p \in (p_2,p_0)$.}
   \vskip -.1cm
 \label{F4m}
\end{figure}

The global $p_2$-branch, which  is
an extension of the positive $p_0$-one
 in Figure \ref{F2m},
 is
shown in Figure \ref{F3m} while the corresponding deformation of
$f$'s in Figure \ref{F4m}. It turns out that it ends up at the
next (even) bifurcation point
 \be
 \label{F3}
  \tex{
  p=p_2= 1+\frac 4{1+2}= \frac 73=2.3333...\, ,
   }
   \ee
so that the branch is expected to be continued for $p<p_2= \frac
73$ in a ``positive" way, etc.

\begin{figure}
\centering
\includegraphics[scale=0.7]{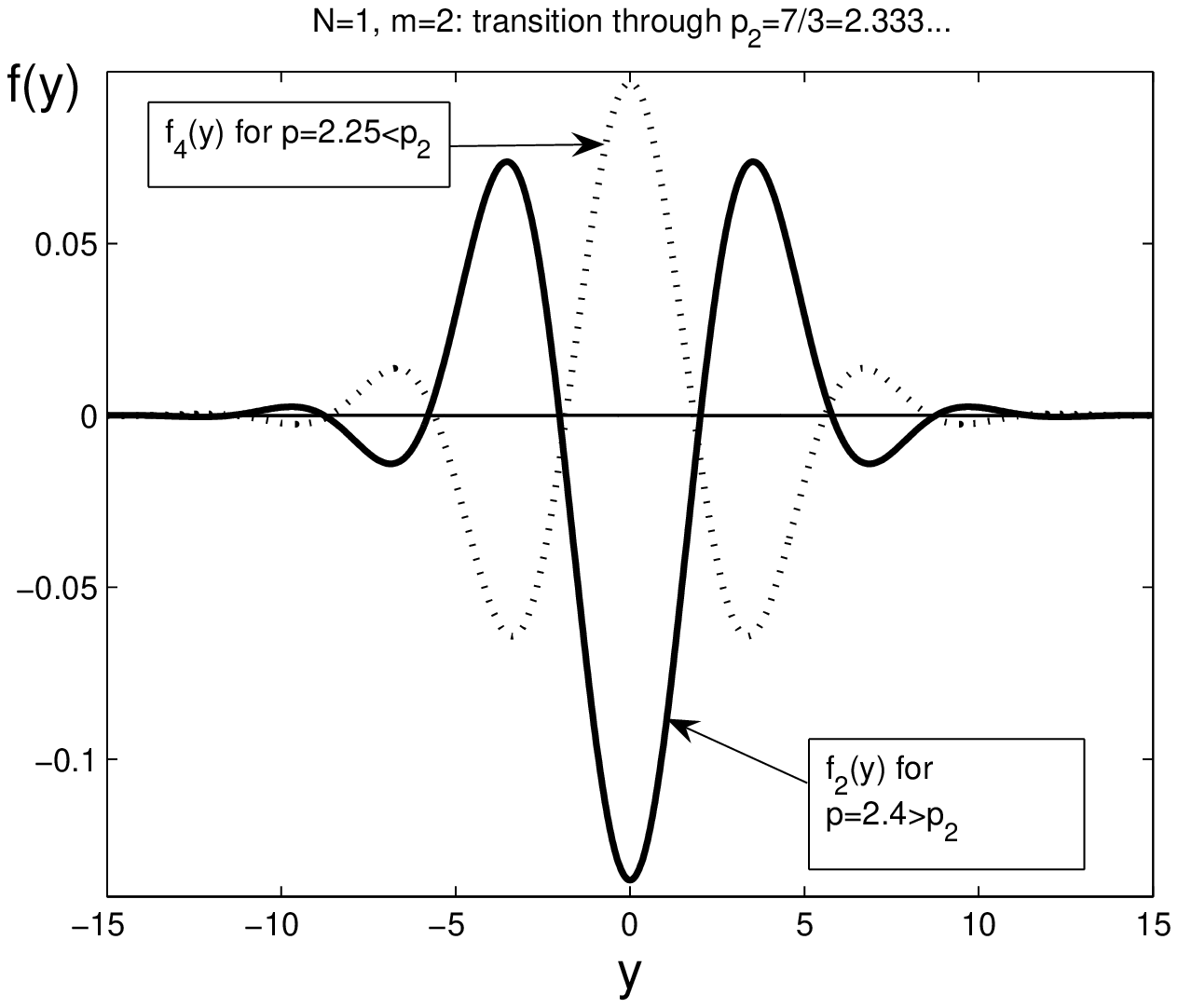}  
\vskip -.5cm \caption{\small Transition of $f_2$ for $p=2.4>p_2$
into $f_4$ for $p=2.25<p_2$, where spatial shape of both mimics
the eigenfunction $\psi_2(y)$ according to \ef{Vexpnm}.}
 \label{F5m}
\end{figure}

\begin{figure}
\centering
\includegraphics[scale=0.7]{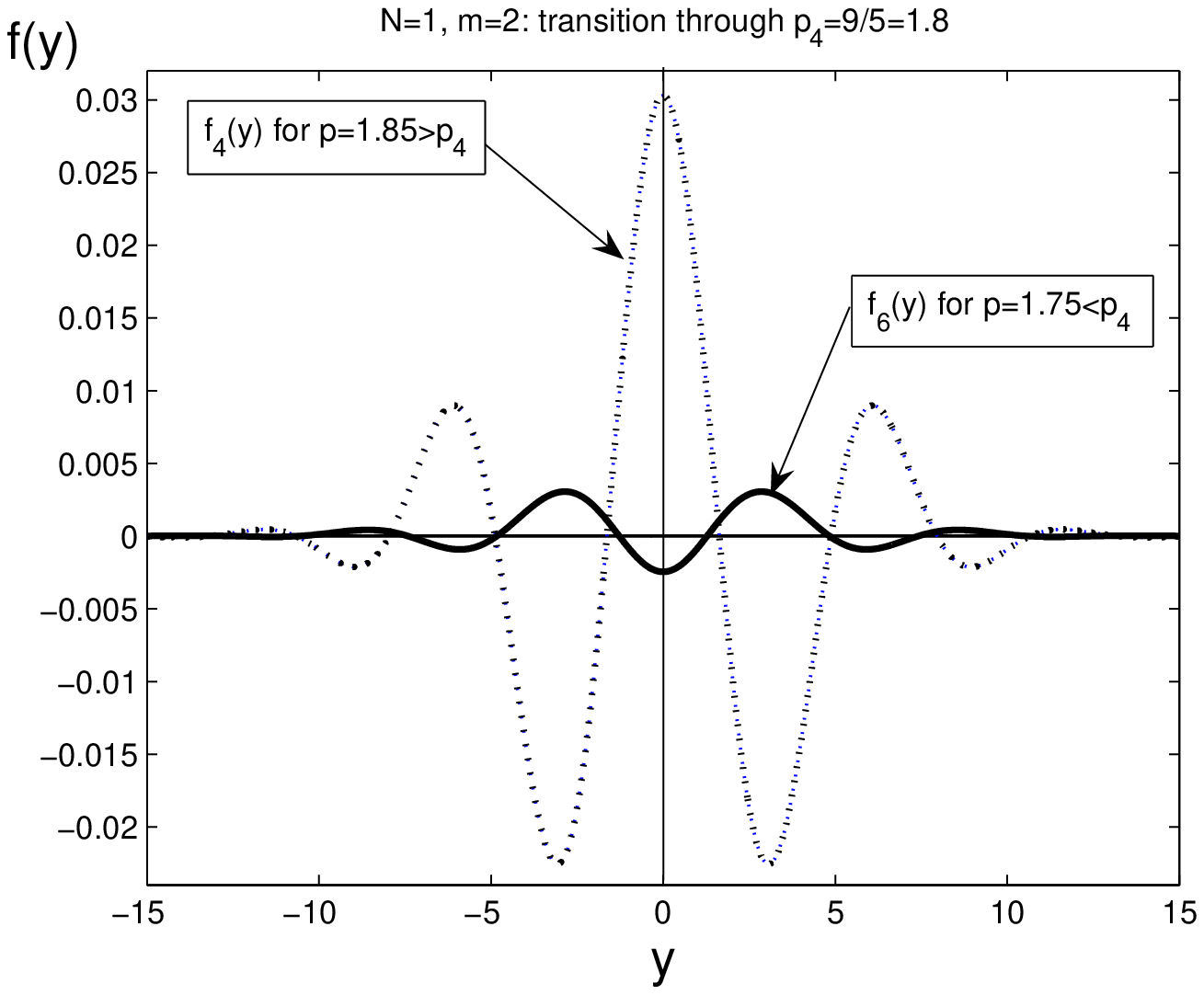}  
\vskip -.5cm \caption{\small Transition of $f_4$ for $p=1.85>p_4$
into $f_6$ for $p=1.75<p_4$, where spatial shape of both mimics
the eigenfunction $\psi_4(y)$ according to \ef{Vexpnm}.}
 \label{Dipp4}
\end{figure}

To justify such transcritical bifurcations at $p=p_l$ for even $l
\ge 2 $, in Figure \ref{F5m}, we present a transition through
$p=p_2= \frac 73\, .$ Similarly, in Figure \ref{Dipp4}, we show
transition from $f_4(y)$ for $p=1.85>p_4=1 + \frac 4{5}=1.8$ to
$f_6(y)$ for $p=1.75<p_4$.



Thus, according to the results given above, we expect that
 there is a continuous deformation along each
connected branch of $f_0$ into $f_2$, $f_2$ into $f_4$, $f_4$ into
$f_6$, etc., i.e., there exists a unique global $p$-branch of even
similarity profiles.
 Hence, we observe that all connected
branches have similar shapes with always two bifurcation points
involved: the right-hand end point $p=p_{4k}$ and the left-hand
end one $p_{4k+2}$.



\subsection {``Approximate" Sturmian zero property}

 Let us comment on the ``Sturmian property" of similarity profiles
$\{f_l(y)\}$. Figures \ref{F11m} and \ref{F5m}  indicate that,
regardless the oscillatory exponential tails, each profile
$f_l(y)$ has a clear ``approximate" (``nonlinear") Sturmian
structure and exhibits  $l+1$ {\em dominant extrema} (meaning $l$
``transversal" zeros in between). In a rigorous mathematical
sense, such properties are known to hold for the second-order
problems. In \cite{GMPOscI}, where \ef{0111} was studied, Sturmian
properties were connected with the category of the functional
subset for each $f_l$, being the corresponding min-max critical
point of the functional, since the category assumes using the
reflection of the functions $(\cdot) \mapsto -(\cdot)$ and hence
the nodal sets of $f_l(y)$ gets more and more complicated as $l$
increases\footnote{There is no still a rigorous treatment of
such zero-set phenomena.}. However, in the present non-variational
case, we cannot use even those rather obscure issues, though
numerical evidence clearly suggests that the approximate Sturmian
nodal properties persist in both variational and non-potential
problems. Honestly, we do not know, which
``mathematical/functional structures" can be responsible for such
a hugely stable Sturmian-like phenomena, and this remains a
challenging open problem of nonlinear operator theory.

 For the higher-order equations, some extra mathematical
reasons for Sturmian properties to persist in an approximate
fashion are discussed in \cite[\S~4.4] {GHUni}. These can be
attributed to the fact that the principal part of (\ref{09Nm})
contains the iteration of two positive operators
 $$
 -D^{
 4}_y= - (-D^2_y)(-D^2_y),
 $$
 and for such pure higher-order operators Sturm's zero property is true
 \cite{Ellias}. Then the linear perturbations affect non-essential
 zeros in the exponential tails only. No rigor justification of such a conclusion
  is available still.


\section{Towards {odd} non-symmetric profiles and their $p$-branches}
 \label{SectN4}

\subsection{An auxiliary discussion}

We begin with noting that, for such fourth-order nonlinear ODEs
\ef{09Nm} with a clearly principal non-coercive operator, for any
$p>1$, we expect to have, at least, a countable set of different
solutions (as in \cite{GHUni, GMPSob, GMPSobII}).

However, in all our previous studies of higher-order elliptic ODE
problems, \cite{GHUni, GMPSob, GMPSobII}, exactly {\bf half} of
such solutions were {\bf odd functions} of $y$ (in 1D; in the
radial symmetry, obviously, odd profiles are not admitted). In the
present case, the profiles $f_1(y)$, $f_3(y)$,...\,,
are not odd (anti-symmetric), since the ODE \ef{09Nm}, unlike that
for \ef{2mGH}, does not admit the corresponding symmetry
 \be
 \label{F2}
f \mapsto -f, \quad y \mapsto -y.
 \ee
 Bearing in mind
 the multiplicity bifurcation results in Section \ref{S2.5}, one
 concludes that other profiles must be non-symmetric in any odd
 sense. Notice that precisely that explains the bifurcation
 expansion \ef{bif12}, where the leading term is odd via
 $\psi_1(y)$, while there exists always an even correction via
 $\psi_2(y)$.


\ssk

We first check some analytical issues concerning such unusual
similarity profiles. Thus, we shoot from $y=+ \iy$ using the same
bundle as \ef{Asyy}, with the coefficients $C^+_{1,2}$,
 \be
 \label{Asyy+}
  \tex{
  f(y) \sim  {\mathrm e}^{-\frac {a_0}2\, y^{4/3}}\big[ C^+_1
   \cos \big(\frac{a_0 \sqrt 3}2 \, y^{\frac 43}\big) + C^+_2
 \sin \big(\frac{a_0 \sqrt 3}2 \, y^{\frac 43}\big)\big] \whereA a_0= 3
 \cdot 2^{-\frac 83}.
  }
  \ee
 Evidently, most of such solutions $f=f(y;C_1^+,C_2^+)$ will
 blow-up at some finite $y_0=y_0(C_1^+,C_2^+)$ according to the following asymptotics:
 as $y \to y_0^+$,
  \be
  \label{Bl1}
   \begin{matrix}
  f^{(4)}=|f|^p(1+o(1)) \LongA
  f(y)=C_0 (y-y_0)^{- \frac 4{p-1}}(1+o(1)),\ssk\\
   \mbox{where} \quad
  C_0^{p-1}= \frac 4{p-1} \big( \frac 4{p-1}+1\big) \big( \frac
  4{p-1}+2\big)\big( \frac 4{p-1}+3\big).
 \end{matrix}
 \ee
Therefore, in order to have a global profile, we have to require
that
 \be
 \label{Bl2}
 y_0(C_1^+,C_2^+)= -\iy.
  \ee

Once we have got such a global solution defined for all $y \in
\re$, we then need to require that, at $y=-\iy$, the algebraic
decay component \ef{Bl4} therein vanishes, i.e.,
 \be
 \label{Bl3}
 C_0^-(C_1^+,C_2^+)=0.
  \ee

We, thus, again arrive at a system of two algebraic equation
\ef{Bl2}, \ef{Bl3}, where the result of Proposition \ref{Pr.Anal}
applies directly to guarantee that {\em the total number of
possible solutions is not more than countable}.

   \ssk



\begin{figure}
\centering
\includegraphics[scale=0.85]{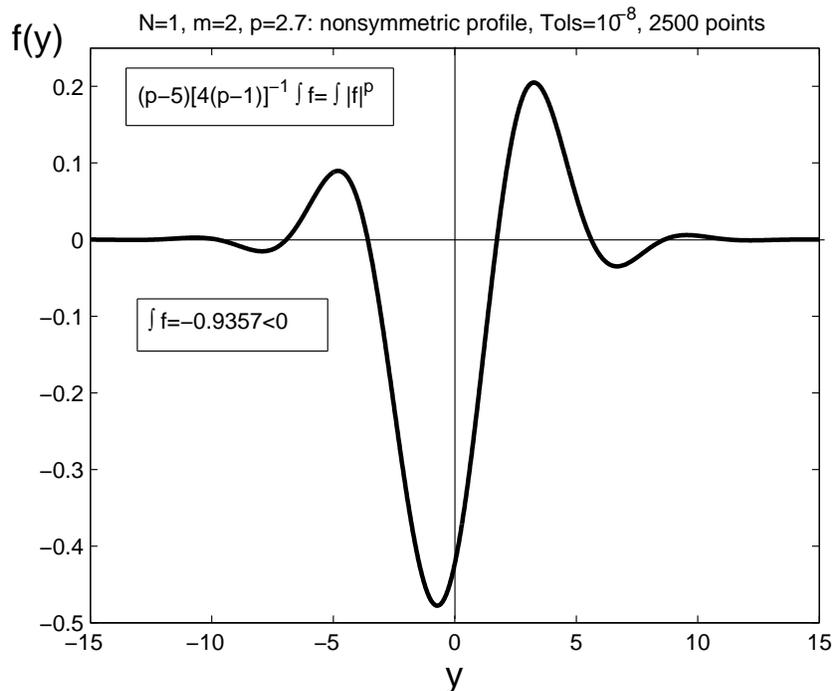}  
\vskip -.5cm \caption{\small A non-symmetric profile $f_3(y)$ for
$p=2.7$,  well satisfying the identity \ef{Int11}.}
   \vskip -.1cm
 \label{D1}
\end{figure}

In Figure \ref{D1}, as a first typical example, we show a
non-symmetric ``dipole-like" profile, denoted by $f_3(y)$ (see
below why such a subscript) for $p=2.7$, which with a sufficient
accuracy $\sim 10^{-5}$ satisfies the identity \ef{Int11}. It
turned out that this identity can be used as a ``blueprint" for
checking the quality in some worse-converging
cases.


\subsection{$p$-branch of $f_3$: from $p_1=3$ to a saddle-node bifurcation}

\begin{figure}
\centering \subfigure[$p$-bifurcation branch of $f_3$]{
\includegraphics[scale=0.5]{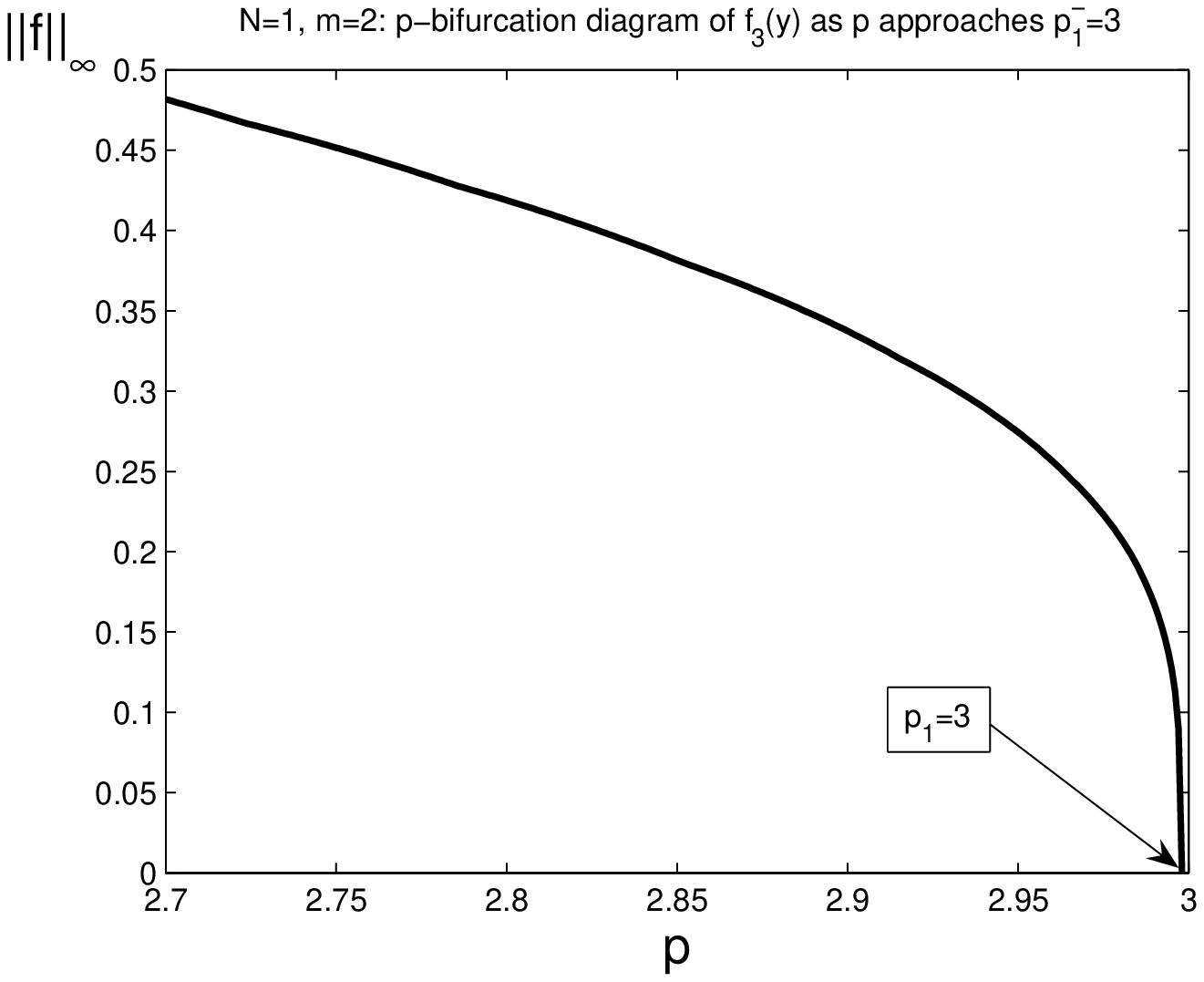} 
} \subfigure[deformation of $f_3$]{
\includegraphics[scale=0.5]{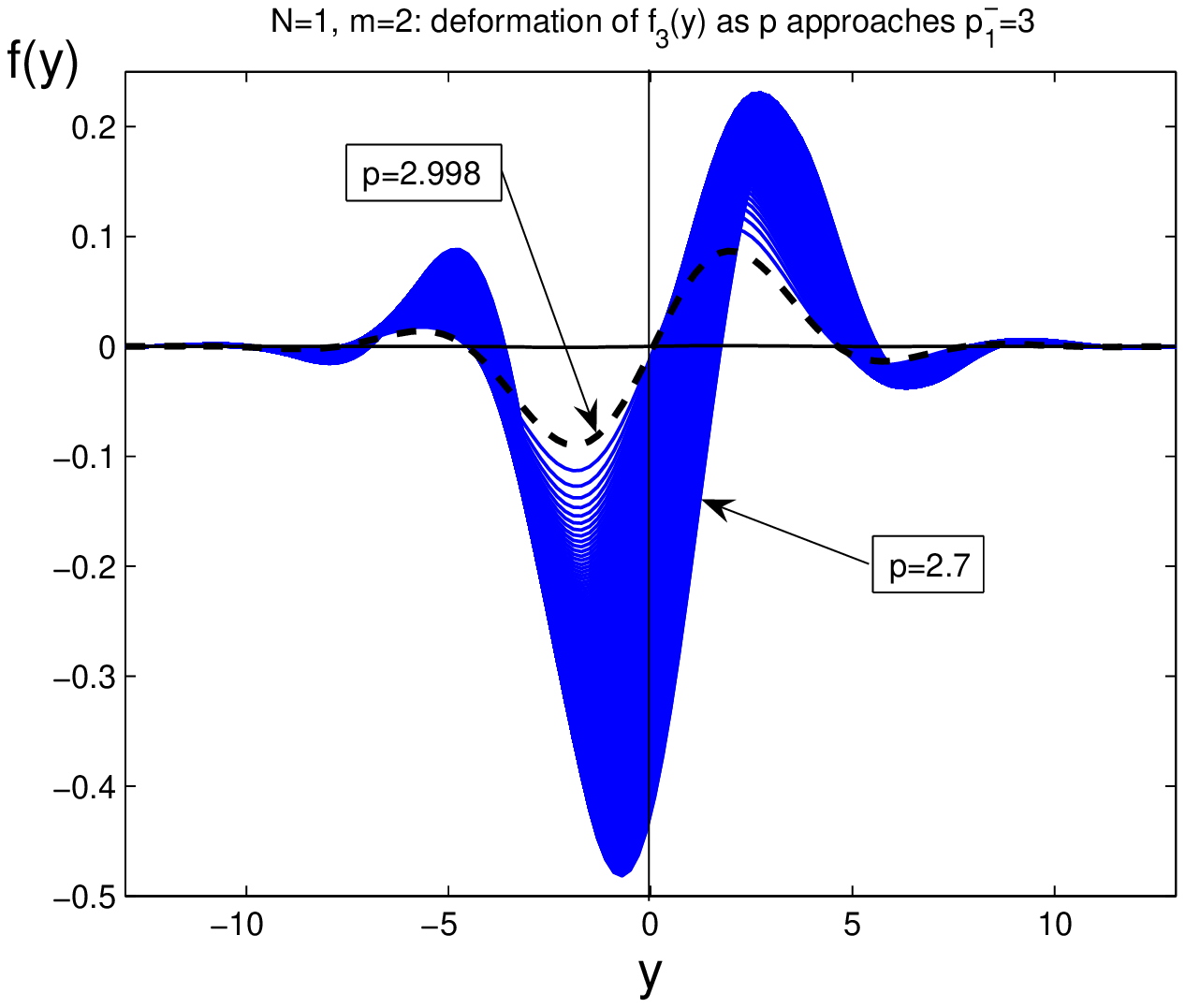}
}
 \vskip -.4cm
\caption{\rm\small $p$-bifurcation branch of $f_3$ (a), and its
deformation (b) for $p \in [2.7,2.924]$.}
  \vskip -.1cm
 \label{DD1}
\end{figure}

Starting from the profile $f_3$ for $p=2.7$  in
 Figure \ref{D1}, we perform a continuation in the parameter $p$
 to get to the bifurcation origin of this $p$-branch.
 Not that surprisingly (cf. Section \ref{S2.5}),
  we observe
in Figure \ref{DD1}(a) that the corresponding bifurcation branch,
with certain accuracy, goes to the odd bifurcation point \ef{h2},
with $l=1$:
 \be
 \label{f31}
  \tex{
  p_1= 1 + \frac 42=3.
  }
  \ee

  In this calculation, we take the continuation step $\D p = 10^{-2}$ (then
  the calculation takes a couple of hours), so, as seen, we cannot approach
  closely to this bifurcation point.  To see approaching $p=p_1$ more clearly,
   we, in addition, took the continuation
  step $\D p= 10^{-4}$ (the calculation of the full branch then took about 16 hours) and observe approaching $p_1=3$ up to
  $p=2.9991$, so that existence of such a ``forbidden" earlier
  bifurcation is without any doubt.
  Figure \ref{DD1}(b) clearly shows that, as $p \to p_1^-=3$, the
  profile $f_3(y)$ (the bold dashed line at $p=2.998$)) takes a typical ``dipole" behaviour governed by
  the second eigenfunction from \ef{eigen}:
   \be
   \label{2nd}
    \tex{
     \psi_1(y)=-  F'(y),
     }
     \ee
where $F(|y|)$ is the even ``bi-harmonic Gaussian" satisfying
\ef{ODEf} for $m=2$. Therefore, this $f_3(y) \sim C(p) \psi_1(y)$,
with some (unknown still) constant $C(p)$, looks like a standard
dipole for the Gaussian for $m=1$
 \be
   \label{2nd1}
    \tex{
     \psi_1(y)=
      \frac 1{2\sqrt{4 \pi}}\, y \, {\mathrm e}^{-y^2/4},
     }
     \ee
but it has oscillatory tails and  further dominant positive and
negative humps.


\begin{figure}
\centering
\includegraphics[scale=0.85]{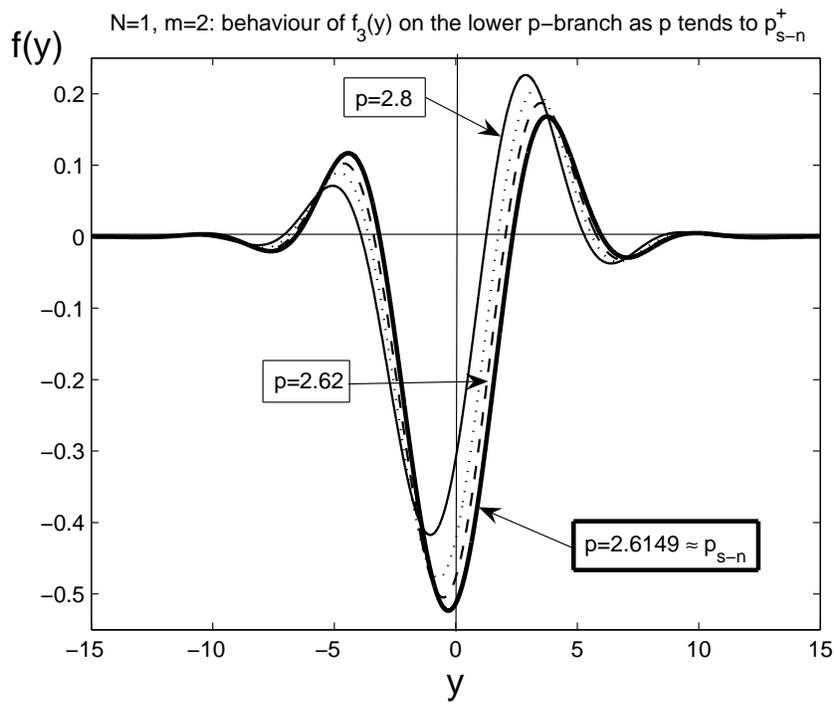}  
\vskip -.5cm \caption{\small Non-symmetric profiles $f_3(y)$ close
to the saddle-node bifurcation \ef{sn1}.}
   \vskip -.1cm
 \label{DD3}
\end{figure}

Extending this $p$-branch of $f_3$ for $p<2.7$, we observe
existence of a saddle-node bifurcation at some $p=p_{s-n}$, where
we obtain the estimate
 \be
 \label{sn1}
  \fbox{$
 2.6148<p_{s-n} \le 2.6149,
 $}
  \ee
by using again the step $\D p=10^{-4}$. The profiles $f_3$ close
to $p_{s-n}$ are shown in Figure \ref{DD3}.

\begin{figure}
\centering
\includegraphics[scale=0.85]{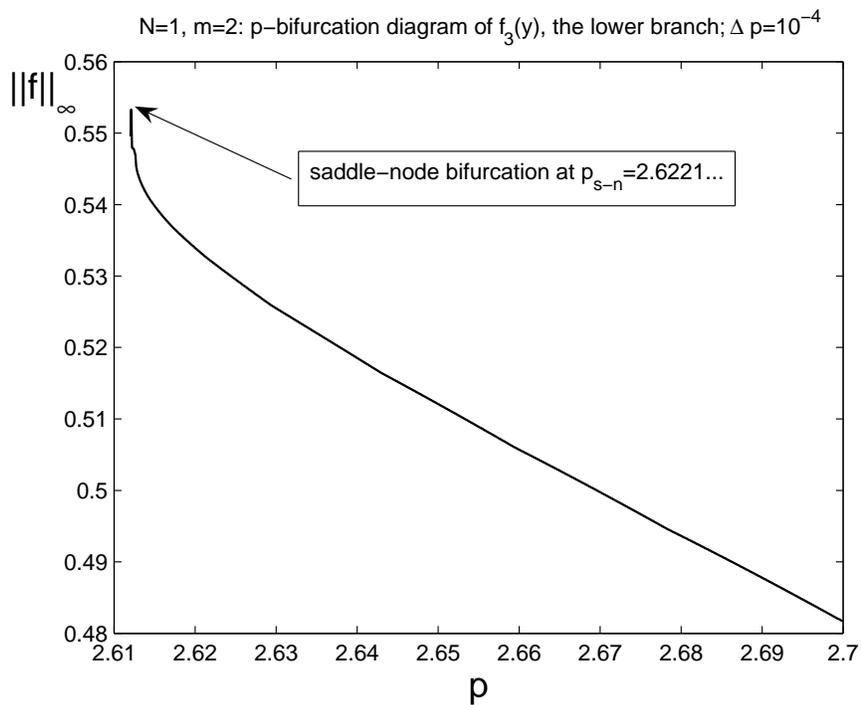}  
\vskip -.5cm \caption{\small The lower bifurcation branch of the
non-symmetric profile $f_3(y)$ for $p \in [2.6149,2.7]$.}
   \vskip -.1cm
 \label{DD4}
\end{figure}

\begin{figure}
\centering
\includegraphics[scale=0.85]{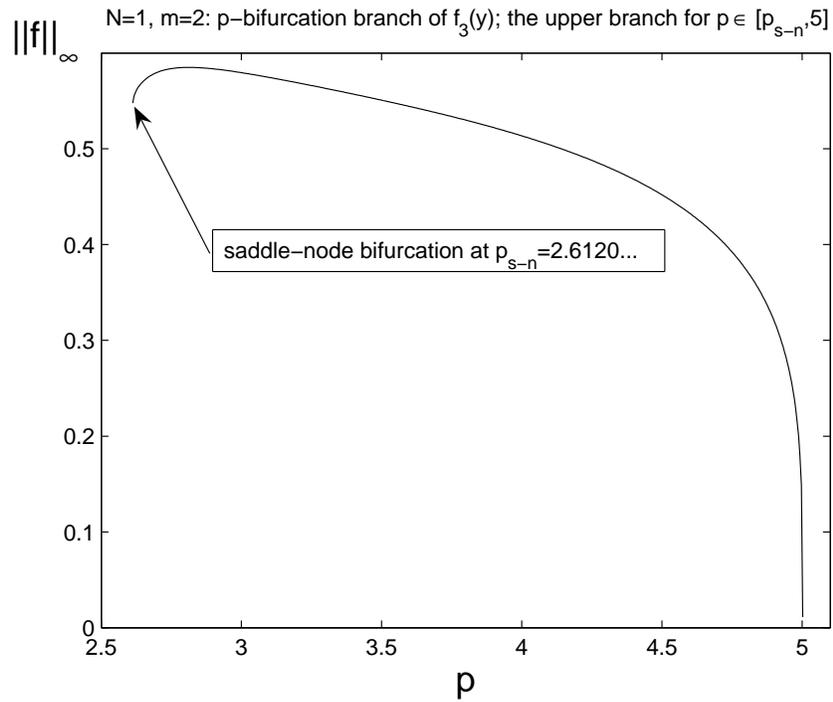}  
\vskip -.5cm \caption{\small The upper bifurcation branch of the
non-symmetric profile $f_3(y)$ for $p \in [2.6149,5]$.}
   \vskip -.1cm
 \label{DD5}
\end{figure}

\begin{figure}
\centering
\includegraphics[scale=0.85]{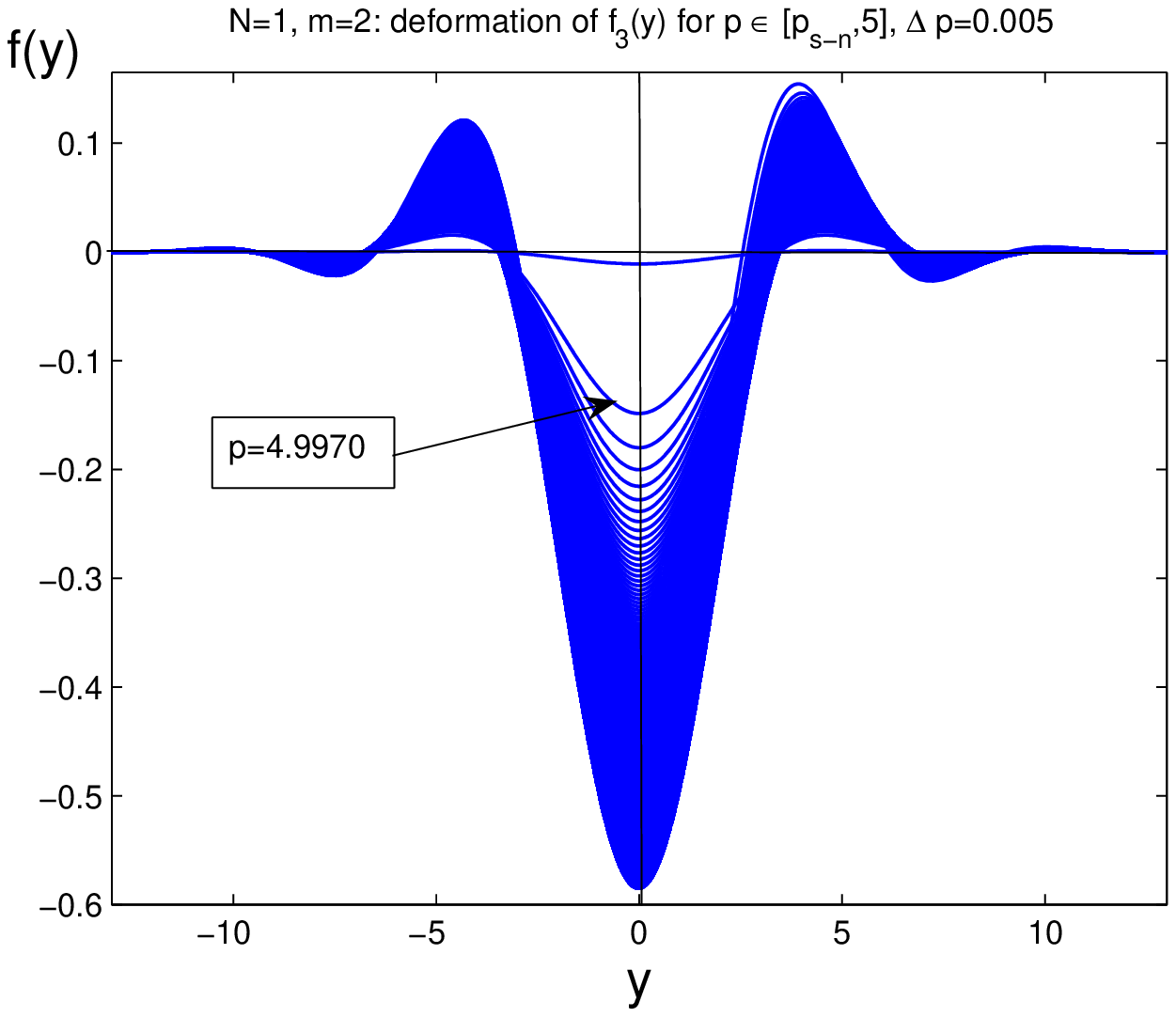}  
\vskip -.5cm \caption{\small Deformation of the non-symmetric
profile $f_3(y)$ for $p \in [2.6149,5]$.}
   \vskip -.1cm
 \label{DD6}
\end{figure}

The {\em lower} bifurcation branch of $f_3$ is shown in Figure
\ref{DD4}. The corresponding {\em upper} bifurcation branch is
shown in Figure \ref{DD5}, while the corresponding deformation of
$f_3$ is presented in Figure \ref{DD6}.

Quite surprisingly, the upper bifurcation branch in Figure
\ref{DD5} ends up at the previous bifurcation point $p_0=5$!
Comparing with Figure \ref{F3m} for symmetric even profiles, we
thus obtain {\sc two} different bifurcation branches (of symmetric
and non-symmetric) solutions originated at $p=5^-$. It is worth
mentioning that the kernel of the linearized operator at $p=5$ is
{\em one-dimensional}. We still do not have a proper explanation
of such a hard and unusual bifurcation phenomenon. However, the
case $l=0$ is not degenerate by \ef{kap00}, so it seems one cannot
create a bifurcation approach similar to that in Section
\ref{S2.5}. This remains an open problems






\section{Second countable family: global linearized patterns}
 \label{S4}




 \subsection{Stable manifold patterns}

  This construction is
 similar to that for $m=1$ and again relies on
stable manifold  \cite{Lun} and Hermitian spectral theory for the
operator pair $\{\BB,\BB^*\}$ in Section \ref{Sect3}. We perform
the same scaling
 of a global solution $u(x,t)$ of \ef{00111} for $t \gg 1$,
 \be
 \label{L1m}
  \tex{
 u(x,t)=t^{-\frac 1{p-1}}\, v(y,\t), \quad y=  x/{t^\frac 1{2m}}, \quad \t= \ln
 t\whereA
  }
  \ee
 \be
 \label{L2m}
  \tex{
  v_\t = {\bf A}(v) \equiv -(-\D)^m v + \frac 1{2m} \, y \cdot \n v + \frac 1{p-1}\, v
  +|v|^{p}, \quad \mbox{so}
  }
   \ee
    \be
    \label{L3m}
     \tex{
    \AAA'(0)= \BB + c_1 I, \quad c_1= \frac 1{p-1}- \frac
    N{2m}=\frac {N(p_0-p)}{2(p-1)} >0 \,\,\, \mbox{for} \,\,\, p<p_0.
     }
     \ee
Thus, $\AAA'(0)$ has the infinite-dimensional stable subspace:
 \be
 \label{L4m}
 E^s= {\rm Span}\{\psi_\b: \,\,\, \l_\b + c_1<0, \,\,
 \mbox{i.e.,} \,\, |\b|>2c_1\}.
  \ee
Using the above spectral properties of $\BB$ \cite{Eg4}, similar
to \cite[\S~5]{GMPOscI}, by invariant manifold theory for
parabolic equations \cite[Ch.~9]{Lun}, we arrive at the following
(see also \cite{GalCr, GHUni}):

\begin{proposition}
 \label{Pr.Linm}

 For any multiindex $\b$ satisfying $|\b|=l>2c_1$, equation $(\ref{L2m})$ admits
 global solutions with the behaviour, as $\t \to +\iy$,
  \be
  \label{L5m}
  v_\b(y,\t)= {\mathrm e}^{(\l_\b+c_1)\t}\varphi_\b(y)(1+o(1))
  \whereA \varphi_\b \in {\rm Span}\{\psi_\b: \,\,\, |\b|=l\},
  \,\,\, \varphi_\b \not = 0.
  \ee
 \end{proposition}

In the original variables \ef{L1m}, the global patterns \ef{L5m}
take the form:
\be
\label{L6m}
 \tex{
 u_\b(x,t) = t^{-\frac{N+|\b|}{2m}} \, \varphi_\b\big( \frac
 x{t^{1/2m}}\big)(1+o(1)) \asA t \to +\iy.
 }
 \ee

  \subsection{Centre manifold patterns}
Unlike the simpler case $m=1$ in \cite[\S~5.2]{GMPOscI}, for the
present  $m \ge 2$, such patterns do exist. Performing a model 1D
analysis of the equation \ef{L2m}, as in
 \cite[\S~5]{GMPOscI}, we conclude that
such patterns may occur if
 \be
 \label{L7}
 \l_\b + c_1=0 \LongA l=|\b|= 2mc_1>0, \quad \mbox{or} \quad p=p_l.
  \ee
  Studying the centre manifold behaviour of the simplest 1D type
   \be
   \label{L8}
   v(\t)= a_l(\t) \psi_l + w^\bot \asA \t \to +\iy,
   \ee
we obtain from \ef{L2m} the following equation for the expansion
coefficient:
 \be
 \label{L9m}
  \dot a_l= \kappa_l \, |a_l|^{p} (1+o(1)) \whereA \kappa_l=
  \langle |\psi_l|^{p}, \psi_l^* \rangle \,\,( \ne 0),
   \ee
which admits global bounded orbits. For instance, noting that
$\kappa_0>0$, one obtains the behaviour
 \be
 \label{SS1}
  a_0(\t)=- [\kappa_0(p-1)]^{-\frac 1{p-1}}\, \t^{-\frac
  1{p-1}}(1+o(1)) \asA \t \to + \iy \quad (p=p_0).
   \ee
Similarly, the same estimate is derived for any $l \ge 0$ provided
that $\kappa_l \not = 0$. Finally, this means that, at such
critical values $p=p_l$, we expect the following logarithmically
perturbed patterns:
 \be
 \label{SS2}
  \tex{
  u_l(x,t) \sim - {\rm sign} \, \kappa_l \, \big[\frac{2m |\kappa_l|}{N+l}\big]^{-\frac
  {N+l}{2m}} \, \big(t \ln t \big)^{-\frac
  {N+l}{2m}} \psi_l\big( \frac x{t^{1/2m}}\big) \asA t \to + \iy.
  }
  \ee

  For the $M$-dimensional eigenspace for $l \ge 1$, we obtain the
  decomposition
   \be
   \label{SS4}
    \tex{
   v(y,\t)= \sum\limits_{|\b|=l} a_\b(\t) \psi_\b(y) + w^\bot(y,\t),
    }
    \ee
that leads to a system of ODEs for the expansion coefficients
$\{a_\b(\t)\}_{|\b|=l}$:
 \be
 \label{SS5}
  \tex{
  \dot a_\g= \big\langle \big|\sum_{|\b|=l} a_\b(\t) \psi_\b \big|^p,
  \psi_\g^* \big\rangle+... \, , \quad |\g|=l.
   }
   \ee
 Assuming the natural ``homogenuity" of this centre subspace
 behaviour:
  \be
  \label{SS6}
   \tex{
   a_\b(\t)= \hat a_\b \, \t^{-\frac {N+l}{2m}}(1+o(1)) \asA \t
   \to +\iy, \quad |\b|=l,
    }
     \ee
     where $\{\hat a_\b\}$ are constants,
\ef{SS5} reduces to an algebraic system (cf. the bifurcation one
\ef{e3m}) of the usual form:
 \be
 \label{SS7}
  \tex{
  \hat a_\g=- \frac{2m}{N+l}  \big\langle \big|\sum_{|\b|=l} \hat a_\b \psi_\b\big|^p,
  \psi_\g^*\big\rangle, \quad |\g|=l.
   }
   \ee
General solvability properties of \ef{SS7}, except some obvious
elementary solutions, and sharp multiplicity results  are unknown.
Of course, as above,  \ef{SS7} is not  variational.

\ssk

\ssk

Again, as for $m=1$,  we arrive at two countable families of
global patterns: the nonlinear  \ef{1.6R} and the linearized
\ef{L6m} ones. Since the rescaled equation \ef{L2m} {\em is not} a
gradient system in any weighted space, their evolution
completeness remains entirely open, though may be expected.






\end{document}